\documentclass[11pt]{article}
\usepackage{amsfonts}
\usepackage{titlesec}
\usepackage{indentfirst}
\usepackage{amsmath}
\usepackage{amscd}
\usepackage{amsthm}
\usepackage{amssymb}
\usepackage{mathrsfs}
\usepackage[dvips]{graphicx}
\usepackage{makeidx}

\newtheorem{thm}{Theorem}

\newtheorem{defn}{Definition}
\newtheorem{conj}{Conjecture}
\newtheorem{remark}{Remark}

\title{New Constructions of  Complex Manifolds  }
\author{ Jinxing Xu \\ \hspace{1em}\textsl{{\small School of Mathematical Sciences, Peking University, Beijing, 100871, P.R. China  }}\\
\hspace{1em}\textsl{{\small E-mail: xujx02@pku.edu.cn}}}

\date{}
\begin {document}

\maketitle
\begin{abstract}
 For a generic anti-canonical hypersurface  in each smooth toric Fano $4-$fold with rank $2$ Picard group , we prove there exist three isolated rational curves
 in it .  Moreover , for all these $4-$folds except one , the contractions of  generic anti-canonical hypersurfaces along the three rational curves can be deformed to smooth threefolds diffeomorphic to connected
 sums of $S^{3}\times S^{3}$ . In this manner, we obtain complex
 structures with trivial canonical bundles on some connected sums of
 $S^{3}\times S^{3}$ . This construction is an analogue of that in Friedman \cite{F2} , Lu and Tian \cite{LT}
 which used only  quintics in $\mathbb{P}^{4}$ .

\end{abstract}

\section{Introduction}
\ This paper is  resulted from an attempt towards the
\emph{Reid's fantasy} for some families of Calabi-Yau threefolds
. Let's first recall some notions .

\begin{defn} [\cite{R}]
Let $Y$ be a Calabi-Yau threefold and $\phi : Y\rightarrow \bar{Y}$
be a birational contraction onto a normal variety . If there exists
a complex deformation (smoothing) of $\bar{Y}$ to a Calabi-Yau
threefold $\tilde{Y}$ , then the process of going from $Y$ to
$\tilde{Y}$ is called a geometric transition and denoted by $T(Y,
\bar{Y} ,\tilde{Y})$ .  A transition $T(Y, \bar{Y} ,\tilde{Y})$ is
called conifold if $\bar{Y}$ admits only ordinary double points as
singularities and the resolution morphism $\phi$ is a small
resolution (i.e. replacing each ordinary double point by a smooth
rational curve).

\end{defn}

Note for a conifold transition $T(Y, \bar{Y} ,\tilde{Y})$ , the
exceptional set of the morphism $\phi$ is  several not intersecting
smooth rational curves each with normal bundle $\mathcal
{O}(-1)\oplus \mathcal {O}(-1)$ in $Y$ and conversely , given some
finite not intersecting smooth rational curves each with normal
bundle $\mathcal {O}(-1)\oplus \mathcal {O}(-1)$ in $Y$ , we can
contract them to get $\bar{Y}$ admitting only ordinary double points
as singularities . The smoothing of $\bar{Y}$ has been studied by
several people . For example , we have the following theorem :

\begin{thm}[Y. Kawamata , G.Tian , see \cite{Ti}]
 \ Let $\bar{Y}$ be a singular threefold with $l$ ordinary double points as
 the only singular points $p_{1} ,\ldots , p_{l}$ . Let $Y $ be a
 small resolution of $\bar{Y}$ by replacing $p_{i}$ by smooth rational curves
 $C_{i}$ . Assume that $Y$ is cohomologically K\"{a}hler and has
 trivial canonical line bundle . Furthermore , we assume that the
 fundamental classes $[C_{i}]$ in $H^{2}(Y ; \Omega^{2}_{Y
 })$ satisfy a relation $\Sigma_{i}\lambda_{i}[C_{i}]=0$ with each
 $\lambda_{i}$ nonzero . Then $\bar{Y}$ can be deformed into a smooth
 threefold $\tilde{Y}$
 .
\end{thm}
 A special case of the above theorem was obtained by R . Friedman in \cite{F1}
 .

 \ The conifold transition process was firstly (locally )observed by
 H .Clemens in \cite{Cle2} , where he explained that locally a
 conifold transition is described by a suitable $S^{3}\times D_{3}$
 to $S^{2}\times D_{4}$ surgery . Roughly speaking , for the conifold
 transition $T(Y, \bar{Y} ,\tilde{Y})$ from $Y$ to $\tilde{Y}$ , it
 kills  $2-$cycles in $Y$ and increases $3-$cycles in $\tilde{Y}$ .
 For a precise relation between their Betti numbers , one can
 consult Theorem $3.2$ in \cite{R} . In  Theorem $1$ ,  if the
 fundamental classes $[C_{i}]$ generates $H^{4}(Y ; \mathbb{C})$ ,
 then we would have $b_{2}(\tilde{Y})=0$ . By results of C.T.C.Wall
 in \cite{Wa} , $\tilde{Y}$ would be diffeomorphic to a connected sum
 of $S^{3}\times S^{3}$ , and
the number of copies is $\frac{b_{3}(\tilde{Y})}{2} + l - b_{2}(Y)$
. We have the following two problems  (cfr. \cite{R}):

\begin{enumerate}
\item  Wether every
projective Calabi-Yau threefold is birational to a Calabi-Yau
threefold $Y$ such that $H^{2}(Y ; \mathbb{C})$ is generated by
rational curves and these curves satisfy the conditions in Theorem
$1$  .\vskip 1pc
\item
Wether the moduli space $\mathcal {N}_{r}$ of complex structures on
the connected sum of $r$ copies of $S^{3}\times S^{3}$ is
irreducible .
\end{enumerate}

\ If we could have positive answers to both of the above problems ,
then we would verify the famous :
\begin{conj}
[the Reid's fantasy  , see \cite{Reid}] Up to some kind of
inductive limit over $r$ , the birational classes of projective
Calabi-Yau threefolds can be fitted together , by means of conifold
transitions , into one irreducible family parameterized by the
moduli space $\mathcal {N}$ of complex structures over suitable
connected sum of copies of $S^{3}\times S^{3}$ .
\end{conj}

In P.S.Green and  T.H\"{u}bsch \cite{GH} , they proved that the
moduli spaces of some Calabi每Yau threefolds , which are complete
intersections in products of projective spaces , were connected each
other by conifold transitions . As for other families of Calabi-Yau
threefolds , for example that which are anti-canonical hypersurfaces
in toric Fano $4-$folds , M. Kreuzer and H. Skarke \cite{KS} proved
they can be connected by geometric transitions , but the transitions
they used are not conifold transitions . So a natural question is
about the Reid's fantasy for these families of Calabi-Yau threefolds
. In this paper , we study the first problem above . Our main result
is :

\begin{thm}
 \ For each  toric smooth Fano $4-$folds $X$  with  rank $2$ Picard group(the only toric Fano $4-$fold with  rank $1$ Picard group is $\mathbb{P}^{4}$ ) , let $Y$ be a generic
 anti-canonical hypersurface of $X$. Then
 there exist three smooth rational curves $C_{i}$ in $Y$ such that
 each one has normal bundle
$\mathcal {O}(-1)\oplus \mathcal {O}(-1)$ in $Y$ ,
 the fundamental classes
$[C_{i}]$ generate $H^{4}(Y;\mathbb{Z})\simeq H_{2}(Y;\mathbb{Z})$
and they satisfy a relation $\Sigma_{i}\lambda_{i}[C_{i}]=0$ with
each $\lambda_{i}$ nonzero . Moreover ,  with the exception of the
variety $B_{1}$ in Batyrev's classification for smooth toric Fano
$4-$folds \cite{Ba}, the three curves do not intersect with each
other .
\end{thm}
Note by Lefschetz's hyperplane theorem , $H_{2}(Y;\mathbb{Z})\simeq
H_{2}(X;\mathbb{Z})$ , so $H_{2}(Y;\mathbb{Z})$ is a free abelian
group with rank $2$ . Since $X$ is simply connected , $H^{2}(X ;
\mathbb{Z})$ is also a free abelian group with rank $2$ . Hence the
classes $[C_{i}]$ generating $H_{2}(Y;\mathbb{Z})$ is equivalent to
generating $H_{2}(X ;\mathbb{Z})$ , and this is also equivalent to
the (dual) cohomology classes represented by $C_{i}$ in $H^{2}(X
;\mathbb{Z})$ generate $H^{2}(X ;\mathbb{Z})$ . Using Theorem $1$
and the discussions above , we can get connected sums of
$S^{3}\times S^{3}$ after contracting these rational curves and
smoothing . The number of copies of $S^{3}\times S^{3}$  is
summarized in the following table , in which the name of toric
$4-$folds is from Batyrev \cite{Ba} . :
\begin{table}[!h]
\tabcolsep 2mm
\begin{center}
\begin{tabular}{r@{}lr@{}lr@{}lr@{}lr@{}lr@{}lr@{}lr@{}lr@{}l}
\hline \multicolumn{2}{c}{Toric Fano $4-$fold
}&\multicolumn{2}{c}{$B_{2}$}&\multicolumn{2}{c}{$B_{3}$}&\multicolumn{2}{c}{$B_{4}$}&\multicolumn{2}{c}{$B_{5}$}&\multicolumn{2}{c}{$C_{1}$}&\multicolumn{2}{c}{$C_{2}$}&\multicolumn{2}{c}{
$C_{3}$}&\multicolumn{2}{c}{$C_{4}$}
\\ \hline
Number of copies of $S^{3}\times S^{3}$&   &104&   &92&   &88&    &88&   &97&   &88&   &88&   &85&  \\

\hline
\end{tabular}
\end{center}
\end{table}

 So according to Theorem $1$ , we get complex structures with trivial
canonical bundles on these connected sums of $S^{3}\times S^{3}$ .
In Lu  and Tian \cite{LT} , they obtained complex structures with
trivial canonical bundles for the connected sum of $m$ copies of
$S^{3}\times S^{3}$ for each $m\geq 2$ , using quintics in
$\mathbb{P}^{4}$ , combing a preceding result of R. Friedman
\cite{F2} .

 The next question is about the relations of the complex structures obtained using
different $4-$folds . For example , in the above table , we observe
that for the varieties $B_{4}$ , $B_{5}$ , $C_{2}$, $C_{3}$ , we can
obtain the same topological connected sum of $S^{3}\times S^{3}$ ,
but we don't know whether the various complex structures on it are
the same , or at least lie in the same deformation class . We also
would like to compare the complex structures obtained using
$\mathbb{P}^{4}$ in Friedman \cite{F2} , Lu and Tian  \cite{LT} and
that obtained using these toric $4-$ folds with rank $2$ Picard
groups . Obviously this question is closely related to the second
problem preceding Conjecture $1$ . Because these complex manifolds
are not K\"{a}hler manifolds (they have vanishing $b_{2}$), very few
techniques are available in dealing with them , and very little is
known about complex structures over connected sums of $S^{3}\times
S^{3}$ . Still and all , the results obtained in this paper
 can be viewed to be a first step for finding connections of
 these Calabi-Yau threefolds and provide  examples of non-intersecting isolated rational curves in some Calabi-Yau threefolds other than quintics (compare T.Johnsen and A.L.Knutsen \cite{JoKn}, T. Johnsen and S. L. Kleiman \cite{JoKl}).

 The paper is organized as follows .

 \ In Sec. $2$ , we  recall the homogenous coordinates on a toric
 variety of D.Cox \cite{Co1}, and give the homogenous coordinates representations for embeddings of $\mathbb{P}^{1}$ to complete nonsingular
toric varieties .

\ In Sec. $3$ , we extend an argument from Clemens \cite{Cle} to
conclude that , if we can find some anti-canonical hypersurface
containing a smooth rational curve with normal bundle $\mathcal
{O}(-1)\oplus \mathcal {O}(-1)$ and with a fixed topological type ,
then generic anti-canonical hypersurfaces will contain a smooth
rational curve with normal bundle $\mathcal {O}(-1)\oplus \mathcal
{O}(-1)$ and with the same  topological type . At the end of this
section , we give a direct method to get the non-intersecting
property for the rational curves in an anti-canonical hypersurface .

\ In Sec. $4$ , we analyze each toric smooth Fano $4-$folds $X$ with
rank $2$ Picard groups  , and for each case we construct three
rational curves $C_{10}$ , $C_{20}$ , $C_{30}$ and three
anti-canonical hypersurfaces $Y_{10}$ , $Y_{20}$ , $Y_{30}$ such
that $C_{i0}$ lies in $Y_{i0}$ with normal bundle $\mathcal
{O}(-1)\oplus \mathcal {O}(-1)$ for $i= 1 , 2 , 3$ . We obtain the
non-intersecting property for the rational curves  , with the
exception of the variety $B_{1}$ . Then using the results in Sec.
$3$ we get our main theorem  .

\textbf{Acknowledgements :} The author would like to sincerely thank
his thesis advisor Professor Gang Tian for proposing this study  and for his continuous encouragement .

\section{Embeddings of $\mathbb{P}^{1}$ to complete nonsingular
toric varieties}

 \ In this section, we will describe all the embeddings of $\mathbb{P}^{1}$
to complete nonsingular toric varieties using the homogeneous
coordinates on a toric variety of D.Cox \cite{Co1}.  First recall
the homogeneous coordinates on a toric variety.

 \ Let $X$ be the toric variety determined by a fan $\Delta$ in $N
\simeq \mathbb{Z}^{n}$. As usual , $M$ will denote the $\mathbb{Z}$
dual of $N$ , and cones in $\Delta$ will be denoted by $\sigma$ .
The one dimensional cones of $\Delta$ form the set $\Delta(1)$ . And
given $\rho\in\Delta(1)$ , let $n_{\rho}$ denote the generator of
$\rho\cap N$ . If $\sigma$ is any cone in $\Delta$ , then
$\sigma(1)=\{\rho\in\Delta(1):\rho\subset\sigma\}$ is the set of one
dimensional faces of the cone $\sigma$ . Throughout of this section
, we assume $X$ to be a complete nonsingular toric variety . So for
each maximal cone $\sigma$ , $\{n_{\rho}: \rho\in\sigma(1)\}$ form a
base of $N$ .

 \ Each $\rho\in\Delta(1)$ corresponds to a $T-$invariant irreducible
Weil divisor $D_{\rho}$ in $X$ . Where $T$ is the torus in $X$ . The
free abelian group of $T-$invariant irreducible Weil divisor on $X$
will be denoted by $\mathbb{Z}^{\Delta(1)}$ .

\ We consider the map

\begin{equation}
M\rightarrow \mathbb{Z}^{\Delta(1)} \texttt{ defined by
}m\rightarrow
D_{m}=\sum_{\rho\in\Delta(1)}<m,n_{\rho}>D_{\rho}\notag
\end{equation}

\ Where $<,>$ means the pairing between elements in $M$ and $N$ .
The map is injective since $\Delta(1)$ spans $N\otimes_{\mathbb{Z}}
\mathbb{R}$ . We have an exact sequence

\begin{equation}
0\rightarrow M\rightarrow \mathbb{Z}^{\Delta(1)} \rightarrow
A_{n-1}(X)\rightarrow 0
\end{equation}

Where $A_{n-1}(X)$ denotes the divisor class group of $X$ . Since
$X$ is complete and nonsingular , $A_{n-1}(X)\simeq Pic(X)$ is a
free abelian group of rank $d-n$ , where $d$ is the number of
one-dimensional faces of $\Delta$ , and $n$ is the dimension of $X$.

For each $\rho\in\Delta(1)$ , we introduce a variable $X_{\rho}$ .
And for each cone $\sigma \in\Delta$ , let
$X_{\hat{\sigma}}=\prod_{\rho \not\in\sigma}X_{\rho}$ , then we can
think of the variety

\begin{equation}
Z=\{x\in \mathbb{C}^{\Delta(1)}: X_{\hat{\sigma}}= 0\texttt{ for all
}\sigma \in \Delta\}\subset \mathbb{C}^{\Delta(1)} \notag
\end{equation}
as the \texttt{"}exceptional\texttt{" }subset of
$\mathbb{C}^{\Delta(1)}$ .

If we apply $Hom( \texttt{-- } ,  \mathbb{C}^{*})$ to the exact
sequence $(1)$ , then we get the exact sequence

\begin{equation}
1\rightarrow G \rightarrow (\mathbb{C}^{*})^{\Delta(1)}\rightarrow T
\rightarrow 1\notag
\end{equation}
where $G=Hom(A_{n-1}(X) , \mathbb{C}^{*})$ and $T=Hom(M ,
\mathbb{C}^{*})$ , both are  products of copies of $\mathbb{C}^{*}$
.

Since $(\mathbb{C}^{*})^{\Delta(1)}$ acts naturally on
$\mathbb{C}^{\Delta(1)}$ , its subgroup
$G\subset(\mathbb{C}^{*})^{\Delta(1)}$ acts on
$\mathbb{C}^{\Delta(1)}$  also . Now $X$ is the geometric quotient
of $\mathbb{C}^{\Delta(1)}-Z$ by $G$ : $X \simeq
\mathbb{C}^{\Delta(1)}-Z \diagup G$

For each $\rho\in\Delta(1)$ , we will canonically  associate to it a
line bundle $L_{\rho}$ and a section $s_{\rho}$ of this line bundle
as follows .

For each maximal cone $\sigma$ of $\Delta$ , denote the
corresponding affine piece of $X$ by $U_{\sigma}\simeq\texttt{ Spec
}\mathbb{C}[\check{\sigma} \cap M]$ , then define $f_{\rho ,
\sigma}\in M $:
\begin{equation}
f_{\rho , \sigma}(n_{\tau})=
\begin{cases}
1& \text{if $\tau=\rho$ },\\
0& \text{otherwise }.
\end{cases}
\end{equation}
for $\tau\in\Delta(1)\cap\sigma$ .

Since $\{n_{\tau}:\tau\in\Delta(1)\cap\sigma\}$ form a base of $N$ ,
$f_{\rho , \sigma}$ is well defined , and $f_{\rho ,
\sigma}\in\check{\sigma}\cap M$ . So $\chi^{f_{\rho , \sigma}}$ is a
regular function on $U_{\sigma}$ . One can check easily that for two
maximal cones $\sigma_{1}$ and $\sigma_{2}$ , $\chi^{f_{\rho ,
\sigma_{1}}}=g_{\sigma_{1}\sigma_{2}}\chi^{f_{\rho , \sigma_{2}}}$
on $U_{\sigma_{1}}\cap U_{\sigma_{2}}$ , where
$g_{\sigma_{1}\sigma_{2}}$ is a nowhere vanishing regular function
on $U_{\sigma_{1}}\cap U_{\sigma_{2}}$ . Using these
$g_{\sigma_{1}\sigma_{2}}$ as transition functions , we ge a line
bundle  $L_{\rho}$ , and those $\chi^{f_{\rho , \sigma}}$ for
maximal cones $\rho$ determine a section $s_{\rho}$  of this line
bundle . It's easy to see that the zero divisor of  $s_{\rho}$ is
just the Weil divisor $D_{\rho}$ .

Generally , if we choose $\rho_{1} ,\rho_{2} ,\ldots ,\rho_{m} \in
\Delta(1)$ , and $d_{1}, d_{2}, \ldots , d_{m} $ integers , we
define
$f_{\sum_{i=1}^{m}d_{i}\rho_{i},\sigma}=\sum_{i=1}^{m}d_{i}f_{\rho_{i},\sigma}$
, and similarly , using
$g_{\sum_{i=1}^{m}d_{i}\rho_{i},\sigma_{1}\sigma_{2}}=\chi^{f_{\sum_{i=1}^{m}d_{i}\rho_{i},\sigma_{1}}}\diagup\chi^{f_{\sum_{i=1}^{m}d_{i}\rho_{i},\sigma_{2}}}$
on the intersection of  $U_{\sigma_{1}}\cap U_{\sigma_{2}}$ for two
maximal cones  $\sigma_{1}$ and $\sigma_{2}$ as transition functions
, we get a meromorphic section $s_{\sum_{i=1}^{m}d_{i}\rho_{i}}$ of
a line bundle $L_{\sum_{i=1}^{m}d_{i}\rho_{i}}$ , which is
isomorphic to $L_{\rho_{1}}^{d_{1}}\otimes
L_{\rho_{2}}^{d_{2}}\otimes\ldots \otimes L_{\rho_{m}}^{d_{m}}$ .
Moreover, this section is regular if $d_{1}, d_{2}, \ldots , d_{m} $
are all nonnegative integers  . one can check that for integers
$d_{1}, d_{2}, \ldots , d_{m} $ and $c_{1}, c_{2}, \ldots , c_{m} $
, if
$\sum_{i=1}^{m}d_{i}D_{\rho_{i}}=\sum_{i=1}^{m}c_{i}D_{\rho_{i}}$ in
$A_{n-1}$ , then
$g_{\sum_{i=1}^{m}d_{i}\rho_{i},\sigma_{1}\sigma_{2}}=g_{\sum_{i=1}^{m}c_{i}\rho_{i},\sigma_{1}\sigma_{2}}$
for any two maximal cones $\sigma_{1}$ and $\sigma_{2}$ . That is ,
the transition functions of $L_{\sum_{i=1}^{m}d_{i}\rho_{i}}$ and
$L_{\sum_{i=1}^{m}c_{i}\rho_{i}}$ coincide  . So in this case , the
quotient $s_{\sum_{i=1}^{m}d_{i}\rho_{i}}\diagup
s_{\sum_{i=1}^{m}c_{i}\rho_{i}}$ is well defined at the points on
which $s_{\sum_{i=1}^{m}c_{i}\rho_{i}}$ is not vanishing .

Using the discussion above , we see that the  sections $s_{\rho}$
for $\rho\in\Delta(1)$ can be used to determine the homogenous
coordinates of $X(\Delta)\simeq \mathbb{C}^{\Delta(1)}-Z \diagup G$

As an application of the homogenous coordinate description for the
complete nonsingular toric variety $X(\Delta)$ , we will obtain all
of the homogenous coordinate representations for the anti-canonical
hypersurfaces on $X$ . Recall the anti-canonical bundle of $X$ is
isomorphic to $L_{\sum_{\Delta(1)}}\rho$ , and a base of  regular
sections is determined by  points in $Q = \{m\in M:<m , n_{\rho}>
\geq -1 , \forall \rho \in \Delta(1)\}$ . Since for any maximal cone
$\sigma$ , the section
$s_{\sum_{\Delta(1)}\rho}|_{U_{\sigma}}=\chi^{\sum_{\rho\in\sigma(1)}f_{\rho,\sigma}}$
, so  when restricted on $U_{\sigma}$ , an anti-canonical  section
represented by $m\in Q$ is equal to
$\chi^{\sum_{\rho\in\sigma(1)}f_{\rho,\sigma}+m}$ . Then using the
homogenous coordinates $(X_{\rho})_{\rho\in\Delta(1)}$ ,
$\chi^{\sum_{\rho\in\sigma(1)}f_{\rho,\sigma}+m}=0$ is equivalent to
$\prod_{\rho\in\Delta(1)}X_{\rho}^{<m,n_{\rho}>+1}=0$ . So any
anti-canonical hypersurface on $X$ has the form $\sum_{m\in
Q}a_{m}\prod_{\rho\in\Delta(1)}X_{\rho}^{<m,n_{\rho}>+1}=0$ , where
$a_{m}\in \mathbb{C}$ are complex numbers .

Now we will describe all the embeddings of $\mathbb{P}^{1}$ to
$X(\Delta)$ . Let $i: \mathbb{P}^{1}\rightarrow X $ be a morphism
from $\mathbb{P}^{1}$ to the toric variety $X$ . Then we have a
homomorphism of their Picard groups : $i^{*}: Pic(X)=A_{n-1}(X)
\rightarrow Pic(\mathbb{P}^{1})=\mathbb{Z}$ . Under this
homomorphism , suppose $i^{*}(L_{\rho})=d_{\rho}$ for $\rho\in
\Delta(1)$ , then the section $s_{\rho}$ is pulled back to a
$d_{\rho}$ form $f_{\rho}(s,t)$ on $\mathbb{P}^{1}$ . So under the
homogenous coordinates on $X$ , the morphism $i$ has the following
form :

\begin{align}
\mathbb{P}^{1}\rightarrow X \notag{}\\
 (s,t)\rightarrow (f_{\rho}(s,t))_{\rho\in\Delta(1)}\notag{}
\end{align}

We call this rational curve  in $X$ has type
$(d_{\rho})_{\rho\in\Delta(1)}$ . It is a generalization of the
concept of degree for rational curves in projective spaces . Note
the homomorphism from $\mathbb{Z}^{\Delta(1)}$ to $\mathbb{Z}$
determined by the integers $d_{\rho} (\rho\in\Delta(1))$ is the
composition of the map $\mathbb{Z}^{\Delta(1)}\rightarrow
A_{n-1}(X)$ in the exact sequence $(1)$ and the homomorphism $i^{*}:
A_{n-1}(X)\rightarrow \mathbb{Z}$ . So  a set of integers
$d_{\rho}(\rho\in\Delta(1))$ is induced by a morphism of
$\mathbb{P}^{1}$ to $X$ if and only if $d_{\rho}(\rho\in\Delta(1))$
satisfy $\sum_{\rho\in\Delta(1)}d_{\rho}n_{\rho}=0$ , where recall
that $n_{\rho}$ is the generator of $\rho\cap N$ .

Next we want to determine the homology class represented by an
embedding of $\mathbb{P}^{1}$ . Since $X$ is complete and
nonsingular , $H^{2}(X, \mathbb{Z})\simeq H_{2n-2}(X ,
\mathbb{Z})\simeq A_{n-1}(X)$ is a finitely generated free abelian
group . So if an embedding of $\mathbb{P}^{1}$ has the form
$f_{\rho}(s,t)(\rho\in\Delta(1))$ , where $f_{\rho}(s,t)$ is a
degree $d_{\rho}$ homogenous form of $s,t$ . Then the cohomology
class in $H^{2}(X, \mathbb{Z})\simeq H_{2n-2}(X , \mathbb{Z})\simeq
A_{n-1}(X)$ represented by this rational curve is
$\sum_{\rho\in\Delta(1)}d_{\rho} [D_{\rho}]$ , where $[D_{\rho}]$ is
the class in $A_{n-1}(X)$ represented by the divisor $D_{\rho}$ .
Note the type and the cohomological class determines each other for
a rational curve in $X$ .

\section{Rational curves in a general anti-canonical hypersurface }

In this section , we  fix a complete nonsingular toric Fano 4-fold
$X=X(\Delta)$ , and use the same notations as the last section . We
will prove that if an anti-canonical hypersurface of $X$ contains a
smooth rational curve $C$ with normal bundle $\mathcal {O}(-1)\oplus
\mathcal {O}(-1)$ , then a generic anti-canonical hypersurface of
$X$ will contain a smooth rational curve with the same type as $C$,
and the normal bundle is also $\mathcal {O}(-1)\oplus \mathcal
{O}(-1)$. More precisely, we have the following theorem .

\begin{thm}
Suppose $Y_{0}$ is an anti-canonical hypersurface in $X$, $C_{0}$ is
a smooth rational curve in $X$ with type
$(d_{\rho})_{\rho\in\Delta(1)}$. Assume $C_{0}$ lies in the smooth
part of $Y_{0}$, and the normal bundle satisfies
$N_{C_{0},Y_{0}}\simeq\mathcal {O}(-1)\oplus \mathcal {O}(-1)$ .Then

\begin{enumerate}
\item
$d_{\rho}\geq -1$ , for any $\rho\in\Delta(1)$ . \vskip 1pc

\item For a generic anti-canonical hypersurface $Y$ (so  $Y$ is
smooth , according to Bertini's theorem ), there is a smooth
rational curve $C$ embedded in $Y$ , such that the type of $C$ in
$X$ is the same as that of $C_{0}$ (so $[C]=[C_{0}]$ in $H_{2}(X ,
\mathbb{Z})$), and the normal bundle satisfies
$N_{C,Y}\simeq\mathcal {O}(-1)\oplus \mathcal {O}(-1)$.
\end{enumerate}
\end{thm}

\begin{proof}
We use an argument analogous to the one used in \cite{Cle}. First of
all we will construct two spaces parameterizing all the
anti-canonical hypersurfaces in $X$ and all the rational curves
embedded in $X$ with type $(d_{\rho})_{\rho\in\Delta(1)}$
respectively . For the anti-canonical hypersurfaces , take $Q$
=$\{w\in M : <w, n_{\rho}>\geq -1, \forall  \rho\in\Delta(1)\}$ . It
is well known that $Q$ is a finite set and we have shown in the last
section that any anti-canonical hypersurface of $X$ has the form
$\sum_{w\in
Q}a_{w}\prod_{\rho\in\Delta(1)}X_{\rho}^{<w,n_{\rho}>+1}=0$, where
$a_{w}\in \mathbb{C}$ are complex numbers , and obviously , not all
the constants $a_{w}$ are zero . Denote $d=\sharp |Q|$ as the number
of elements in $Q$ . Then we can take $\mathbb{P}^{d-1}$ as a
parameter space for all the anti-canonical hypersurfaces in $X$ .

For all the rational curves embedded in $X$ with type
$(d_{\rho})_{\rho\in\Delta(1)}$ , note under the homogenous
coordinates of $X$ , any such rational curve has the form
$(f_{d_{\rho}}(s,t))_{\rho\in\Delta(1)}$ , where $s,t$ are the
homogenous coordinates of $\mathbb{P}^{1}$ , and $f_{d_{\rho}}(s,t)$
is a homogenous polynomial of $s,t$ with degree $d_{\rho}$ . By
convention , $f_{d_{\rho}}(s,t)\equiv 0$ if $d_{\rho}<0$ . Let
$\Delta_{*}(1)=\{\rho\in \Delta(1): d_{\rho}< 0\}$,
$\Delta^{*}(1)=\{\rho\in \Delta(1): d_{\rho}\geq 0\}$. Suppose
$f_{d_{\rho}}= \sum_{i=0}^{d_{\rho}}b_{\rho, i}s^{i}t^{d_{\rho}-i}$
, for $\rho \in \Delta^{*}(1)$ . Then it is natural to collect all
the coefficients $b_{\rho, i}$ for $\rho \in \Delta^{*}(1)$, $i=0,
\ldots , d_{\rho}$ to construct a parameterizing space for all the
rational curves with type $(d_{\rho})_{\rho\in\Delta(1)}$ in $X$.
Next we give the precise definition . Recall if
$A_{n-1}(X)\simeq\mathbb{Z}^{m}$, then $X\simeq
\mathbb{C}^{\Delta(1)}-Z\diagup (\mathbb{C}^{*})^{m}$ , where
$Z=\{x\in \mathbb{C}^{\Delta(1)}:
X_{\hat{\sigma}}=\prod_{\rho\not\in \sigma}X_{\rho}= 0\texttt{ , for
all }\sigma \in \Delta\}\subset \mathbb{C}^{\Delta(1)}$ , and
$(\mathbb{C}^{*})^{m}$ acts on $\mathbb{C}^{\Delta(1)}$ in the form
$(\lambda_{1},\ldots, \lambda_{m})\cdot
(X_{\rho})_{\rho\in\Delta(1)}=(\varphi_{\rho}(\lambda_{1},\ldots,\lambda_{m})X_{\rho})_{\rho\in\Delta(1)}$,
 with $\varphi_{\rho}:(\mathbb{C}^{*})^{m}\rightarrow
 \mathbb{C}^{*}$ a homomorphism for  $\rho\in\Delta(1)$. Now
 define

 \begin{equation}
 \mathcal{M}^{'}=
 \mathbb{C}^{\sharp(\Delta_{*}(1))}\times\mathbb{C}^{\sum_{\rho\in\Delta^{*}(1)}(d_{\rho}+1)}-Z^{'}
 \diagup(\mathbb{C}^{*})^{m}\notag
\end{equation}
 where
 \begin{equation}
\begin{split}
Z^{'}=\{(b_{\rho,0})_{\rho\in\Delta_{*}(1)}\times(b_{\rho,0},\ldots,
b_{\rho,d_{\rho}})_{\rho\in\Delta^{*}(1)}: \prod_{\rho\not\in
\sigma}b_{\rho,i(\rho)}= 0, \forall \sigma \in \Delta, \forall
i(\rho) \texttt{ such that } \\
0\leq i(\rho)\leq d_{\rho}\texttt{ if }\rho\in\Delta^{*}(1) \texttt{
, and }i(\rho)=0 \texttt{ if }\rho\in\Delta_{*}(1)\} \subseteq
\mathbb{C}^{\sharp(\Delta_{*}(1))} \times
\mathbb{C}^{\sum_{\rho\in\Delta^{*}(1)}(d_{\rho}+1)}
\end{split}\notag
\end{equation}
 , and
$(\mathbb{C}^{*})^{m}$ acts on
$\mathbb{C}^{\sharp(\Delta_{*}(1))}\times\mathbb{C}^{\sum_{\rho\in\Delta^{*}(1)}(d_{\rho}+1)}$
in the form
\begin{equation}
\begin{split}
(\lambda_{1},\ldots, \lambda_{m})\cdot
(b_{\rho,0})_{\rho\in\Delta_{*}(1)}\times(b_{\rho,0},\ldots,
b_{\rho,d_{\rho}})_{\rho\in\Delta^{*}(1)}= \\
(\varphi_{\rho}(\lambda_{1},\ldots,
\lambda_{m})b_{\rho,0})_{\rho\in\Delta_{*}(1)}\times(\varphi_{\rho}(\lambda_{1},\ldots,
\lambda_{m})b_{\rho,0},\ldots, \varphi_{\rho}(\lambda_{1},\ldots,
\lambda_{m})b_{\rho,d_{\rho}})_{\rho\in\Delta^{*}(1)}
\end{split}\notag
\end{equation}
Now define $\mathcal {M}$ to be the subvariety of $\mathcal {M}^{'}$
with $b_{\rho,0}=0$ for all $\rho\in\Delta_{*}(1)$. It's not hard to
verify that $\mathcal {M}^{'}$ is a nonsingular complete toric
variety with dimension $4+\sum_{\rho\in\Delta^{*}(1)}d_{\rho}$ , and
$\mathcal {M}$ is a nonsingular subvariety of $\mathcal {M}^{'}$
with dimension
$4+\sum_{\rho\in\Delta^{*}(1)}d_{\rho}-\sharp(\Delta_{*}(1))$.

Consider the incidence variety
\begin{equation}
I=\{(a,b)\in \mathbb{P}^{d-1}\times\mathcal {M}:
F_{a}(f^{b}_{d_{\rho}}(s,t))_{\rho\in\Delta(1)}\equiv0\}\subseteq
\mathbb{P}^{d-1}\times\mathcal {M}\notag
\end{equation}
where for $a=(a_{w})_{w\in Q}\in\mathbb{P}^{d-1}$ ,
$F_{a}(X_{\rho})_{\rho\in\Delta(1)}=\sum_{w\in
Q}a_{w}\prod_{\rho\in\Delta(1)}X_{\rho}^{<w,n_{\rho}>+1}$. And for
$b=(0)_{\rho\in\Delta_{*}(1)}\times(b_{\rho,0},\ldots,
b_{\rho,d_{\rho}})_{\rho\in\Delta^{*}(1)}\in\mathcal {M}$,
$f^{b}_{d_{\rho}}(s,t))\equiv0$ if $\rho\in\Delta_{*}(1)$ , and
$f^{b}_{d_{\rho}}(s,t))=\sum_{i=0}^{d_{\rho}}b_{\rho,
i}s^{i}t^{d_{\rho}-i}$ , for $\rho \in \Delta^{*}(1)$.

Using the equality $\sum_{\rho\in\Delta(1)}d_{\rho}<w, n_{\rho}>=0$
for any $w\in M$ , we see that
$F_{a}(f^{b}_{d_{\rho}}(s,t))_{\rho\in\Delta(1)}$ is a homogenous
polynomial of $s,t $ with degree $\sum_{\rho\in\Delta(1)}d_{\rho}$ ,
if  it is  not $0$  . Hence elementary dimension theory implies that
every irreducible component of $I$ has dimension not less than
\begin{equation}
 dim\mathbb{P}^{d-1}
+dim \mathcal {M} -1- \sum_{\rho\in\Delta(1)}d_{\rho}
=dim\mathbb{P}^{d-1} +
4+\sum_{\rho\in\Delta^{*}(1)}d_{\rho}-\sharp(\Delta_{*}(1)) -1-
\sum_{\rho\in\Delta(1)}d_{\rho}
\end{equation}

\ On the other hand , the existence of $C_{0}$ and $Y_{0}$ in the
hypothesis implies there is a point $(a_{0}, b_{0})\in I$ , where
$a_{0}$ is the coefficients of the defining equation of $Y_{0}$, and
$b_{0}$ denotes a   parameterization for $C_{0}$ .  Since the normal
bundle of $C_{0}$ in $Y_{0}$ has no nonzero sections, $C_{0}$ is
infinitesimally rigid in $Y_{0}$ . This implies that the fibre
dimension of the projection $I \rightarrow \mathbb{P}^{d-1}$ at
$(a_{0}, b_{0})$ is exactly $3$ . In fact the fiber is parametrized
by $PGL(2)$. So taking an irreducible component $I_{0}$ of $I$ going
through $(a_{0}, b_{0})$ , we have
\begin{equation}
dimI_{0}\leq dim \mathbb{P}^{d-1} + 3
\end{equation}
Now $(3)$and $(4)$ implies
\begin{equation}
dim\mathbb{P}^{d-1} +
4+\sum_{\rho\in\Delta^{*}(1)}d_{\rho}-\sharp(\Delta_{*}(1)) -
\sum_{\rho\in\Delta(1)}d_{\rho}-1\leq dimI_{0}\leq dim
\mathbb{P}^{d-1} + 3
\end{equation}
Since $\sum_{\rho\in\Delta^{*}(1)}d_{\rho} +
\sum_{\rho\in\Delta_{*}(1)}d_{\rho}=\sum_{\rho\in\Delta(1)}d_{\rho}$
, and $d_{\rho}\geq 0$ for $\rho\in\Delta^{*}(1)$ , $d_{\rho}\leq
-1$ for $\rho\in\Delta_{*}(1)$ . We conclude that all the
inequalities in $(5)$ are in fact equalities , in particular
$d_{\rho}\geq -1$ , for any $\rho\in\Delta(1)$ . That proves the
first claim of the theorem .

Now $dimI_{0}= dim \mathbb{P}^{d-1} + 3$ and that the fibre
dimension of the composed morphism $I_{0}\hookrightarrow I
\rightarrow \mathbb{P}^{d-1}$ at $(a_{0}, b_{0})$ is exactly $3$
will imply that this morphism $I_{0}\rightarrow \mathbb{P}^{d-1}$ is
surjective , and the generic fiber has dimension $3$ . This will
imply that for generic hypersurface $Y$ , there is a rational curve
$C$ embedded in $Y$ , such that the type of $C$ in $X$ is the same
as that of $C_{0}$ , and the normal bundle satisfies
$N_{C,Y}\simeq\mathcal {O}(-1)\oplus \mathcal {O}(-1)$. The
smoothness of $C$ comes from $C_{0}$ is smooth and that to be a
regular embedding is an open condition on $\mathcal {M}$. So we have
proven the second claim in the theorem .

\end{proof}

\begin{remark}
Using some basic deformation theory , one can prove a similar result
replacing the toric variety $X$ by any complete smooth Fano $4-$fold
.
\end{remark}

\begin{remark}
In the above theorem, suppose there is a Zariski open subset $U$ of
$\mathcal {M}$ such that $C_{0}$ lies in $U$(more precisely, there
is a parametrization of $C_{0}$ in  $\mathcal {M}$ which lies in
$U$). Then the rational curve $C$ in the above theorem can be chosen
to lie in $U$ , too . This can be easily seen in the proof.
\end{remark}

According to the above theorem , for a nonsingular complete Fano
toric variety $X$ , if we can find some smooth rational curves
$C_{10}$, $\ldots$, $C_{l0}$ in $X$ such that their types are all
different with each other , and  for each curve $C_{i0}$,$1\leq
i\leq l$ , there is an anti-canonical hypersurface $Y_{i0}$ of $X$
going through $C_{i0}$ such that they satisfy the hypothesis of the
above theorem. Then  for generic anti-canonical hypersurface $Y$ .
$Y$ contains smooth rational curves $C_{1}$, $\ldots$, $C_{l}$ such
that they all have normal bundles $\mathcal {O}(-1)\oplus \mathcal
{O}(-1)$ in $Y$ , and the type of $C_{i}$ is equal to that of
$C_{i0}$ ( so the cohomological class represented by $C_{i}$ is
equal to that represented by $C_{i0}$ ), for $1\leq i\leq l$. So if
the  cohomological classes $[C_{10}]$, $\ldots$, $[C_{l0}]$ satisfy
the condition in  Theorem $1$ , then the same is true for the curves
$C_{i}$  .

By  Remark $2$ , if we choose Zariski open subset $U_{i}$ $(1\leq
i\leq l)$ of $\mathcal {M}$ such that $C_{i0}$ lies in $U_{i}$ for
each $1\leq i\leq l $, then the curves $C_{1}$, $\ldots$, $C_{l}$ in
$Y$ can be also chosen to lie in $U_{1}$, $\ldots$, $U_{l}$
respectively. In practice, we usually choose $U_{i}$ to be the
Zariski subset of $\mathcal {M}$ which represents all the regular
embeddings of $\mathbb{P}^{1}$, or at the same time, some homogenous
coordinates of $X$ have no zero points on the embedded rational
curve.

Using the same notation as the last paragraph , note in the
hypothesis of Theorem $1$ , we require the smooth rational curves
$C_{1}$, $\ldots$, $C_{l}$ lying in $Y$ do not intersect each other.
So next we want to analyze when can we guarantee that for generic
anti-canonical hypersurface $Y$, the  rational curves $C_{1}$,
$\ldots$, $C_{l}$ lying in $Y$ do not intersect each other. Take any
two of these curves , suppose they are $C_{1}$, $C_{2}$ without loss
of generality. We fix two Zariski open set $U_{1}$, $U_{2}$ of the
corresponding parametrizing space $\mathcal {M}_{1}$, $\mathcal
{M}_{2}$ such that $C_{i}$ lies in $U_{i}$ and every point in
$U_{i}$ represents a regular embedding, for $i=1,2$. Consider the
subvariety $U_{12}$ of $\mathcal {M}_{1}\times\mathcal {M}_{2}$:
\begin{equation}
U_{12}=\{(b_{1}, b_{2})\in U_{1}\times U_{2} :C_{1}\cap C_{2}\neq
\emptyset\}\notag
\end{equation}
Where $C_{1}, C_{2}$ denotes the rational curves represented by
$b_{1},b_{2}$ respectively .

Consider the following incident variety :
\begin{equation}
\begin{split}
J= \{(a, b_{1}, b_{2})\in \mathbb{P}^{d-1}\times\ U_{12}:
F_{a}(f^{b_{1}}_{d_{\rho}}(s,t))_{\rho\in\Delta(1)}=\\
F_{a}(f^{b_{2}}_{d_{\rho}}(s,t))_{\rho\in\Delta(1)}\equiv0\}
\subseteq \mathbb{P}^{d-1}\times\mathcal {M}_{1}\times\mathcal
{M}_{2} \notag
\end{split}
\end{equation}

Roughly speaking, $J$ represents the configuration that two
intersecting rational curves lying in an anti-canonical
hypersurface.

Now all we want to do is to find conditions to guarantee the
dimension of the image of  $Pr_{1} $ is strictly less than  $dim$
$\mathbb{P}^{d-1}$ (where $Pr_{1} $ denotes the natural projection
morphism form $J$ to $\mathbb{P}^{d-1}$), for then the image of
$Pr_{1} $ is a lower dimensional constructible set (i.e. finite
union of locally closed set) of $\mathbb{P}^{d-1}$, so the closure
of this image is a lower dimensional closed subvariety in
$\mathbb{P}^{d-1}$ . Note each fibre of the natural projection
morphism form $J$ to $\mathbb{P}^{d-1}$ has dimension not less than
$6$ , because of  the free action of $PGL(2)$ on each of the two
rational curves. So it suffices to prove $dim$ $  J < $ $dim$
$\mathbb{P}^{d-1} + 6$. Considering the natural projection morphism
$Pr_{2}$ from $J$ to $U_{12}$ , we only have to prove that each
fibre of $Pr_{2}$ has dimension strictly less than $dim$ $
\mathbb{P}^{d-1} + 6 - dim$ $U_{12}$ . By the definition of $J$, for
each point $(b_{1}, b_{2})\in U_{12}$, the  fibre  of $Pr_{2}$ at
$(b_{1}, b_{2})$ is the  linear subspace of $\mathbb{P}^{d-1}$ such
that its points represent exactly the anti-canonical hypersurfaces
containing both of  the rational curves represented by $b_{1}$ and
$b_{2}$. So if for any rational curve pair $(C_{1}, C_{2})$
represented by a point in $U_{12}$, we can  find $dim$ $U_{12} - 5 $
anti-canonical sections of $X$ , such that  the hypersurface
corresponding to any nonzero linear combination of these
anti-canonical sections never contain $C_{1}$ and $C_{2}$ at the
same time, then we would be done . That is what we will do for each
concrete toric variety in the next section. Unfortunately , this
method fails in the case of the variety $B_{1}$ . That's why in our
main theorem we have an exceptional case .

\section{Examinations for toric Fano $4-$folds with  rank $2$ Picard groups }

In this section , we give homogenous coordinates representations for
each toric Fano $4-$folds with Picard group rank $2$ . Using these
representations , we will find three smooth rational curves
$C_{10},C_{20},C_{30}$ and anti-canonical hypersurfaces
$Y_{10},Y_{20},Y_{30}$ such that $[C_{i0}]$ satisfy the conditions
in Theorem $1$ , and $C_{i0}$ is embedded in the smooth part of
$Y_{i0}$ with normal bundle $\mathcal {O}(-1)\oplus \mathcal
{O}(-1)$ . Then Theorem $3$ implies a generic hypersurface $Y$ will
contain three smooth rational curves $C_{1},C_{2},C_{3}$ such that
$[C_{i}]$ satisfy the conditions in Theorem $1$ , and $C_{i}$ is
embedded in  $Y$ with normal bundle $\mathcal {O}(-1)\oplus \mathcal
{O}(-1)$ for $i= 1  , 2 , 3$ . Similar to the definition of $U_{12}$
at the end in the last section , we can define $U_{ij}$
 parameterizing intersecting rational curves $C_{i} , C_{j}$ and
analyze the dimension of the space $U_{ij}$ . Finally we prove that
a generic anti-canonical hypersurface does not contain intersecting
rational curves with our chosen topological types .

\subsection{The toric variety $B_{1}$}

\ In the classification of Batyrev \cite{Ba} , the toric variety
$B_{1}$ is defined by a fan $\Delta$ in $\mathbb{R}^{4}$ such that
$\Delta(1) = \{ v_{1} , \ldots , v_{6}\}$ , $\Delta(1)$ generates
$\mathbb{Z}^{4}$ , and elements in $\Delta(1)$ satisfy the following
linear relations (cf. \cite{Ba}) :
\begin{equation}
\left\{ \begin{aligned}
         \ v_{1}+v_{2}+v_{3}+v_{4}=3 v_{6} \\
                  \ v_{5}+v_{6}=0
                          \end{aligned} \right.\notag
\end{equation}

$\mathbb{C}^{*}\times \mathbb{C}^{*}$ acts on this toric variety as
follows :

\begin{equation}
(X_{1}, \ldots , X_{6} )\rightarrow (\lambda_{1}X_{1},\ldots ,
\lambda_{1}X_{4}, \lambda_{1}^{3}\lambda_{2}X_{5},
\lambda_{2}X_{6})\notag
\end{equation}
For $(\lambda_{1}, \lambda_{2})\in \mathbb{C}^{*}\times
\mathbb{C}^{*}$ .

 Under this action , $X\simeq
\mathbb{C}^{6}-Z\diagup \mathbb{C}^{*}\times \mathbb{C}^{*}$ , where
$Z=\{X_{1}X_{5}=X_{2}X_{5}=X_{3}X_{5}=X_{4}X_{5}
=X_{1}X_{6}=X_{2}X_{6}=X_{3}X_{6}=X_{4}X_{6}=0\}$

The anti-canonical forms of $X$ are linear combinations of the
following forms :

\begin{equation}
X_{5}^{2}f_{1}(X_{1},\ldots , X_{4}) , X_{5}X_{6}f_{4}(X_{1},\ldots
, X_{4}),  X_{6}^{2}f_{7}(X_{1},\ldots , X_{4}). \notag
\end{equation}
where $f_{i}(X_{1},\ldots , X_{4})$ denotes a  degree $i$ homogenous
form of $X_{1},X_{2} ,X_{3}, X_{4}$ , for $i\geq 1$.

 Now $H^{2}(X, \mathbb{Z})\simeq H_{6}(X ,
\mathbb{Z})\simeq A_{3}(X)$ is a rank $2$ free abelian group , and
$[D_{4}],[D_{5}]$ form a base of this group .

Consider the rational curve  $(0,0,0,1,s,t)$ in $X$ , which has the
type $(0,0,0,0,1,1)$ . Its cohomology class is $[D_{5}]$ ,  and this
rational curve is embedded in the anti-canonical hypersurface $
X_{1}X_{5}^{2}+X_{2}X_{4}^{3}X_{5}X_{6}+X_{3}X_{4}^{6}X_{6}^{2}=0$
with normal bundle $\mathcal {O}(-1)\oplus \mathcal {O}(-1) $ .

We explain the computation of the normal bundle .

Denote $Y$ as the anti-canonical hypersurface $
X_{1}X_{5}^{2}+X_{2}X_{4}^{3}X_{5}X_{6}+X_{3}X_{4}^{6}X_{6}^{2}=0$ ,
 and $i : \mathbb{P}^{1}\rightarrow Y$ the embedding of the rational
curve $(0,0,0,1,s,t)$ . Then we have two exact sequences :

\begin{equation}
N_{Y, X}^{*} \rightarrow \Omega_{X\diagup \mathbb{C}}\rightarrow
\Omega_{Y\diagup \mathbb{C}}\rightarrow 0
\end{equation}
Where $N_{Y,  X}^{*}$ is the conormal bundle of the hypersurface $Y$
in the toric $4-$fold $X$ . Pull back this exact sequence to
$\mathbb{P}^{1}$ using $i$ , we get :
\begin{equation}
i^{*}(N_{Y, X}^{*}) \rightarrow i^{*}(\Omega_{X\diagup
\mathbb{C}})\rightarrow i^{*}(\Omega_{Y\diagup
\mathbb{C}})\rightarrow 0
\end{equation}

Since the rational curve lies in the smooth part of $Y$ , we have
the exact sequence :

\begin{equation}
0\rightarrow N_{\mathbb{P}^{1}, Y}^{*} \rightarrow
i^{*}(\Omega_{Y\diagup \mathbb{C}})\rightarrow
\Omega_{\mathbb{P}^{1}\diagup \mathbb{C}}\rightarrow 0
\end{equation}
Where $N_{\mathbb{P}^{1}, Y}^{*}\simeq Hom (N_{\mathbb{P}^{1}, Y} ,
\mathcal {O}_{\mathbb{P}^{1}})$ is the conormal bundle of the
rational curve in $Y$ .

Using the two exact sequence $(7)$ and $(8)$ , we can compute
concretely the rank 2 locally free sheaf $N_{\mathbb{P}^{1}, Y}
\simeq \mathcal {O}(-1)\oplus \mathcal {O}(-1)$

Similarly , the rational curve $(s ,t,0,0,0,1)$ has the type
$(1,1,1,1,3,0)$ . Its cohomology class is $[D_{4}]+3[D_{5}]$ , and
this rational curve is embedded in the anti-canonical section
$X_{3}X_{1}^{6}X_{6}^{2}+
X_{4}X_{1}^{4}X_{2}^{2}X_{6}^{2}+X_{5}X_{2}^{4}X_{6}=0$ with normal
bundle $\mathcal {O}(-1)\oplus \mathcal {O}(-1) $ .

We summarize the results in the following table .

\begin{table}[!h]
\tabcolsep 0.5mm \caption{$B_{1}$}
\begin{center}
\begin{tabular}{r@{}lr@{}lr@{}lr@{}l}
\hline \multicolumn{2}{c}{rational curve
$C_{i0}$}&\multicolumn{2}{c}{type}&\multicolumn{2}{c}{cohomology
class}&\multicolumn{2}{c}{anti-canonical hypersurface $Y_{i0}$}
\\ \hline
(0,0,0,1,s,t&)   &(0,0,0,0,1,1&)&[$D_{5}$&]&$X_{1}X_{5}^{2}+X_{2}X_{4}^{3}X_{5}X_{6}+X_{3}X_{4}^{6}X_{6}^{2}=$&0 \\
(s,t,0,0,0,1&)&(1,1,1,1,3,0&)&[$D_{4}$]+3[$D_{5}$&]&$X_{3}X_{1}^{6}X_{6}^{2}+
X_{4}X_{1}^{4}X_{2}^{2}X_{6}^{2}+X_{5}X_{2}^{4}X_{6}=$&0\\
(s,t,0,0,$s^{4}$,t&)&(1,1,1,1,4,1&) &[$D_{4}$]+4[$D_{5}$&]  & $X_{3}X_{5}^{2}+X_{4}X_{1}^{6}X_{6}^{2}+X_{5}X_{6}X_{2}^{4}-X_{6}^{2}X_{2}^{3}X_{1}^{4}=$&0\\
\hline
\end{tabular}
\end{center}
\end{table}
Where in each row , the computation shows $C_{i0}$ is embedded in
the smooth part of  $Y_{i0}$ with normal bundle $\mathcal
{O}(-1)\oplus \mathcal {O}(-1)$ .

Now the three curves $(0,0,0,1,s,t)$ , $(s,t,0,0,0,1)$,
$(s,t,0,0,s^{4},t)$ are denoted by $C_{10}$ ,$C_{20}$, $C_{30}$
respectively , and the corresponding parameterizing spaces are
denoted by  $\mathcal {M}_{1}$,$\mathcal {M}_{2}$,$\mathcal {M}_{3}$
. Take the Zariski open set $U_{1}$ of $\mathcal {M}_{1}$ such that
points in $U_{1}$ represent regular embedded rational curves on
which the homogenous $X_{4}$ does not have zero points . So $U_{1}$
parameterizes exactly smooth rational curves in $X$ having the form
$(c_{1},c_{2},c_{3},c_{4},\alpha_{1}(s,t),\beta_{1}(s,t))$ , where
$\alpha_{1}(s,t),\beta_{1}(s,t)$ denote linear forms of $s,t$ ,
$c_{i} (1\leq i\leq4)$ are constants and $c_{4}$ is not zero.
Similarly, take $U_{2}$ as the open set of $\mathcal {M}_{2}$
parameterizing regular embedded rational curves on which the
homogenous coordinate $X_{6} $ does not have zero points . Take
$U_{3}$ to be the open set of $\mathcal {M}_{3}$ parameterizing the
regular embedded rational curves . Using the similar notation as the
end of last section, we want to  analyze $dim$ $
U_{ij}(i,j=1,2,3,i\neq j)$ and find $dim $ $U_{ij}-5 $
anti-canonical sections such that any nonzero linear combinations of
these sections will not contain both of the  rational curves
represented by points in $U_{ij}$. This will suffice to prove that
generic anti-canonical hypersurface will contain smooth rational
curves in $U_{1}$,$U_{2}$,$U_{3}$ , and these curves do not
intersect each other. We will use $C_{1}$,$C_{2}$,$C_{3}$ to denote
curves represented by points in $U_{1}$,$U_{2}$,$U_{3}$
respectively.

The case of $C_{1}$ and $C_{2}$

By definition, $C_{1}$ has a representation with the form $(c_{1},
c_{2},c_{3}, 1, \alpha_{1}(s,t),\beta_{1}(s,t))$, where $c_{i}(1\leq
i \leq 3)$ are constants , and $\alpha_{1}(s,t),\beta_{1}(s,t)$ are
degree $1$ homogenous forms of $s,t$. $C_{2}$ has a representation
with the form
$(\alpha_{1}(s,t),\beta_{1}(s,t),\gamma_{1}(s,t),\delta_{1}(s,t),\alpha_{3}(s,t),1)$
, where the subscript numbers denote the degrees of the
corresponding homogenous forms of $s,t$. Now consider $dim $
$U_{12}$. In the representation of $C_{2}$, there are $12$
coefficients in the homogenous forms
$\alpha_{1}(s,t),\beta_{1}(s,t),\gamma_{1}(s,t),\delta_{1}(s,t),\alpha_{3}(s,t)$
, and modulo the action of $\mathbb{C}^{*}\times \mathbb{C}^{*}$ ,
the appearance of $C_{2}$ will contribute $11$ to $dim$ $U_{12}$.
When fixing $C_{2}$, since $C_{1}$ has to intersect with $C_{2}$,
$(c_{1}, c_{2},c_{3}, 1)$ has to lie in the rational curve
$(\alpha_{1}(s,t),\beta_{1}(s,t),\gamma_{1}(s,t),\delta_{1}(s,t))$,
this contributes $1$ to $dim$ $U_{12}$. At last, the $4$
coefficients of $\alpha_{1}(s,t),\beta_{1}(s,t)$ in the
representation of $C_{1}$ modulo the action of $\mathbb{C}^{*}$ will
contribute $3$ to $dim$ $U_{12}$ . So we get $dim$ $U_{12}=15$, and
we have to find $10$ anti-canonical sections satisfying the
condition we just required.

Note that an invertible linear substitution of the homogenous
coordinates $X_{1}$ , $X_{2}$ ,$X_{3}$ ,$X_{4}$ induces an
automorphism of $X$, so for any point $(b_{1}, b_{2})\in U_{12}$ ,
the two rational curves $C_{1}$, $C_{2}$ represented by $b_{1} ,
b_{2}$ can be assumed to be $(0,0,0,1,s,t)$ and
$(0,0,s,t,\alpha_{3}(s,t),1)$ , after choosing appropriate
homogenous coordinates on $\mathbb{P}^{1}$. Then its easy to verify
the following $10$ anti-canonical forms satisfy our requirement.
\begin{equation}
X_{5}^{2}X_{4},X_{5}X_{6}X_{4}^{4},X_{6}^{2}X_{4}^{7},X_{6}^{2}X_{3}X_{3}^{i}X_{4}^{6-i}\
(i=0,\ldots, 6) . \notag
\end{equation}

The case of $C_{1}$ and $C_{3}$

Similar to the analysis in the last paragraph, $dim $ $U_{13}=17$ ,
and without loss of generality, we can assume $C_{1}=(0,0,0,1,s,t),
C_{3}=(0,0,s,t,\alpha_{4}(s,t),\alpha_{1}(s,t))$, since $C_{3} $ is
a smooth rational curve , the homogenous forms
$\alpha_{4}(s,t),\alpha_{1}(s,t)$ have no common factors in
$\mathbb{C}[s,t]$. Pick a degree $3$ homogenous form
$\beta_{3}(s,t)$ such that it has no common factors with
$\alpha_{1}(s,t)$, then it's easy to verify that  the following $12$
anti-canonical forms will suffice :
\begin{equation}
X_{5}^{2}X_{4}, X_{5}X_{6}X_{4}^{4}, X_{6}^{2}X_{4}^{7},
X_{5}^{2}X_{3}, X_{5}X_{6}X_{3}\beta_{3}(X_{3}, X_{4}),
X_{6}^{2}X_{3}X_{3}^{i}X_{4}^{6-i}(i=0, \ldots, 6) .\notag
\end{equation}

The case of $C_{2}$ and $C_{3}$

In this case , our method fails . So we have an exceptional case in
the main theorem .
\subsection{The toric variety $B_{2}$}

This variety  is defined by a fan $\Delta$ such that elements in
$\Delta(1)$ satisfy :
\begin{equation}
\left\{ \begin{aligned}
         \ v_{1}+v_{2}+v_{3}+v_{4}=2v_{6} \\
                  \ v_{5}+v_{6}=0
                          \end{aligned} \right. \notag
\end{equation}

\ $\mathbb{C}^{*}\times \mathbb{C}^{*}$ acts on this toric variety
as follows :

\begin{equation}
(X_{1}, \ldots , X_{6} )\rightarrow (\lambda_{1}X_{1},\ldots ,
\lambda_{1}X_{4}, \lambda_{1}^{2}\lambda_{2}X_{5},
\lambda_{2}X_{6})\notag
\end{equation}

Under this action , $X\simeq \mathbb{C}^{6}-Z\diagup
\mathbb{C}^{*}\times \mathbb{C}^{*}$ , where
$Z=\{X_{1}X_{5}=X_{2}X_{5}=X_{3}X_{5}=X_{4}X_{5}
=X_{1}X_{6}=X_{2}X_{6}=X_{3}X_{6}=X_{4}X_{6}=0\}$

The anti-canonical forms of $X$ are linear combinations of the
following forms :

\begin{equation}
X_{5}^{2}f_{2}(X_{1},\ldots , X_{4}) , X_{5}X_{6}f_{4}(X_{1},\ldots
, X_{4}),  X_{6}^{2}f_{6}(X_{1},\ldots , X_{4}). \notag
\end{equation}
where $f_{i}(X_{1},\ldots , X_{4})$ are degree $i$ homogenous forms
of $X_{1},X_{2} ,X_{3}, X_{4}$ , for $i\geq 1$.

 Now $H^{2}(X, \mathbb{Z})\simeq H_{6}(X ,
\mathbb{Z})\simeq A_{3}(X)$ is a rank $2$ free abelian group , and
$[D_{4}],[D_{5}]$ form a base of this group .

\begin{table}[!h]
\tabcolsep 0.1mm \caption{$B_{2}$}
\begin{center}
\begin{tabular}{r@{}lr@{}lr@{}lr@{}l}
\hline \multicolumn{2}{c}{rational curve
$C_{i0}$}&\multicolumn{2}{c}{type}&\multicolumn{2}{c}{cohomology
class}&\multicolumn{2}{c}{anti-canonical section $Y_{i0}$}
\\ \hline
(0,0,0,1,s,t&)&(0,0,0,0,1,1&)&[$D_{5}$&]&$X_{1}X_{5}^{2}X_{4}+X_{2}X_{4}^{3}X_{5}X_{6}+X_{3}X_{4}^{5}X_{6}^{2}=$&0 \\
(s,t,0,0,0,1&)&(1,1,1,1,2,0&)&[$D_{4}$]+2[$D_{5}$&]&$X_{3}X_{1}^{5}X_{6}^{2}+
X_{4}X_{1}^{3}X_{2}^{2}X_{6}^{2}+X_{5}X_{2}^{4}X_{6}=$&0\\
(s,t,0,0,$s^{3}$,t&)&(1,1,1,1,3,1&) &[$D_{4}$]+3[$D_{5}$&]  & $X_{3}X_{1}X_{5}^{2}+X_{4}X_{2}^{5}X_{6}^{2}+X_{5}X_{6}X_{1}^{2}X_{2}^{2}-X_{6}^{2}X_{2}X_{1}^{5}=$&0\\
\hline
\end{tabular}
\end{center}
\end{table}

Denote the three rational curves in the first column by $C_{10}$ ,
$C_{20}$ , $C_{30}$ respectively , and denote the corresponding
parameterizing spaces by $\mathcal {M}_{1}$ , $\mathcal {M}_{2}$ ,
$\mathcal {M}_{3}$ , then take the corresponding three Zariski open
sets $U_{1}$ , $U_{2}$ , $U_{3}$ such that the rational curves
parameterized by points in $U_{i}$ are smooth . We will use the same
notation as the last case and proceed in a similar way.

The case of $C_{1}$ and $C_{2}$

$dim$ $U_{12}=14$ , $C_{1}, C_{2}$ can be assumed to be
$(0,0,0,1,s,t) , (0,0,s,t,\alpha_{2}(s,t),1)$ , and the following
$9$ anti-canonical forms will suffice :
\begin{equation}
X_{5}^{2}X_{4}^{2}, X_{5}X_{6}X_{4}^{4},
X_{6}^{2}X_{3}^{i}X_{4}^{6-i}(i=0, \ldots, 6) .\notag
\end{equation}

The case of $C_{1}$ and $C_{3}$

$dim$ $U_{13}=16$ , $C_{1} , C_{3}$ can be assumed to be
$(0,0,0,1,s,t) , (0,0,s,t,\alpha_{3}(s,t),\alpha_{1}(s,t))$ , and
the following $11$ anti-canonical forms will suffice :
\begin{equation}
X_{5}^{2}X_{4}^{2}, X_{5}X_{6}X_{4}^{4}, X_{6}^{2}X_{4}^{6},
X_{5}^{2}X_{3}\hat{\beta}_{1}(X_{3},X_{4}),
X_{5}X_{6}X_{3}\hat{\beta}_{3}(X_{3},X_{4}),X_{6}^{2}X_{3}X_{3}^{i}X_{4}^{5-i}(i=0,
\ldots, 5) . \notag
\end{equation}
Where $\hat{\beta}_{1}(s,t), \hat{\beta}_{3}(s,t)$ are homogenous
forms with degree $1, 3$ respectively and neither one has  common
factors with $\alpha_{1}(s,t)$ . This can be used to guarantee that
when restricting on the rational curve $C_{3}$ , the $8$ forms
\begin{equation}
X_{5}^{2}X_{3}\hat{\beta}_{1}(X_{3},X_{4}),
X_{5}X_{6}X_{3}\hat{\beta}_{3}(X_{3},X_{4}),X_{6}^{2}X_{3}X_{3}^{i}X_{4}^{5-i}(i=0,
\ldots, 5)\notag
 \end{equation}
are linearly independent as degree $8$ homogenous forms of $s ,t$ .

The case of $C_{2}$ and $C_{3}$

$dim$ $U_{23}=20$ , the rational curve pair $C_{2} , C_{3}$ can be
assumed to be :
\begin{equation}
C_{2}=(\alpha_{1}(s,t),\beta_{1}(s,t),\gamma_{1}(s,t),\delta_{1}(s,t),\alpha_{2}(s,t),1)
\notag
\end{equation}
\begin{equation}
C_{3}=(\tilde{\alpha}_{1}(s,t),\tilde{\beta}_{1}(s,t),\tilde{\gamma}_{1}(s,t),\tilde{\delta}_{1}(s,t),\alpha_{3}(s,t),\tau_{1}(s,t))
\notag
\end{equation}

Note
$(\alpha_{1}(s,t),\beta_{1}(s,t),\gamma_{1}(s,t),\delta_{1}(s,t))$ ,
$(\tilde{\alpha}_{1}(s,t),\tilde{\beta}_{1}(s,t),\tilde{\gamma}_{1}(s,t),\tilde{\delta}_{1}(s,t))$
represent two rational curves in $\mathbb{P}^{3}$. Denote them by
$C_{2}^{'}$ , $C_{3}^{'}$ respectively . Recall the notation at the
end of Section $2$ , if for the morphism form $J$ to $U_{23}$ ,
there is a rational curve pair $C_{2}$ , $ C_{3}$ in the image of
$J$ such that $C_{2}^{'}$ does not coincide with $ C_{3}^{'}$ as
lines in $\mathbb{P}^{3}$ , then modulo a re-parameterization and an
automorphism of $X$ , $C_{2}$ , $C_{3}$ can be assumed to be :
\begin{equation}
C_{2}=(0,s,0,t,\alpha_{2}(s,t),1),
C_{3}=(0,0,s,t,\alpha_{3}(s,t),\alpha_{1}(s,t)) .\notag
\end{equation}
And the following $15$ anti-canonical forms will suffice :
\begin{equation}
 X_{6}^{2}X_{2}^{i}X_{4}^{6-i}(i=0,
\ldots, 6), X_{5}^{2}X_{3}\hat{\beta_{1}}(X_{3},X_{4}),
X_{5}X_{6}X_{3}\hat{\beta}_{3}(X_{3},X_{4}),X_{6}^{2}X_{3}X_{3}^{i}X_{4}^{5-i}(i=0,
\ldots, 5) .\notag
\end{equation}
Where the homogenous forms  $\hat{\beta}_{1}$ , $\hat{\beta}_{3}$
are the same
 as those in the case of $C_{1}$ and $C_{3}$ . If for any rational curve pair $C_{2}, C_{3}$ in the image of $J$ ,
$C_{2}^{'}$ coincides with $ C_{3}^{'}$ as lines in $\mathbb{P}^{3}$
, then the dimension of this image $dim $ $Im(J)\leq 18$ , and
modulo a re-parameterization and an automorphism of $X$ , $C_{2}$ ,
$C_{3}$ can be assumed to be :
\begin{equation}
C_{2}=(0,0,s,t,\alpha_{2}(s,t),1),
C_{3}=(0,0,s,t,\alpha_{3}(s,t),\alpha_{1}(s,t)) .\notag
\end{equation}
Since $\alpha_{1}(s,t)$ has no common factors with $\alpha_{3}(s,t)$
, we can find degree $3$ homogenous forms $\gamma_{3}(s,t) ,
\delta_{3}(s,t)$ such that the $4$ homogenous forms
$\alpha_{1}\alpha_{2}$ , $\alpha_{3}$, $\gamma_{3}$, $\delta_{3}$
are linearly independent as homogenous forms of $s,t$ . Now consider
the $7$ anti-canonical forms :

\begin{equation}
\begin{split}
X_{5}^{2}\alpha_{1}^{2}(X_{3},X_{4}),X_{5}X_{6}\alpha_{1}(X_{3},X_{4})\alpha_{3}(X_{3},X_{4}),
X_{5}X_{6}\alpha_{1}(X_{3},X_{4})\gamma_{3}(X_{3},X_{4}), \\
X_{5}X_{6}\alpha_{1}(X_{3},X_{4})\delta_{3}(X_{3},X_{4}),X_{6}^{2}\alpha_{3}(X_{3},X_{4})\alpha_{3}(X_{3},X_{4}),\\
X_{6}^{2}\alpha_{3}(X_{3},X_{4})\gamma_{3}(X_{3},X_{4}),X_{6}^{2}\alpha_{3}(X_{3},X_{4})\delta_{3}(X_{3},X_{4}).
\end{split}
\end{equation}
When restricted to $C_{3}$, they reduce to $3$ homogenous forms of
$s,t$ with degree $8$ . Since the anti-canonical forms could
generate all the degree $8$ forms on $C_{3}$ , we can pick $6$
anti-canonical forms $f_{i} (1\leq i \leq 6 )$ such that with  the
$7$ forms in $(9)$ , these $13$ anti-canonical forms generate all
the degree $8$ forms on $C_{3}$ when restricted on it . Next it's
direct to verify that when restricted on $C_{2}$ , the $7$ forms in
$(11)$ generate a dimension $7$, $6$, or $5$ linear subspace of
degree $6$ forms on $C_{2}$ , depending on $\alpha_{2}(s,t)$ and
$\alpha_{3}(s,t)$ has a degree $0$ , $1$ , or $2$ greatest common
divisor respectively . According to the above analysis , we can
always find $dim$ $ Im(J) - 5$ anti-canonical forms such that any
nonzero linear combination of them does not contain $C_{2}$ and
$C_{3}$ at the same time , for $C_{2}$ and $C_{3}$ represented by
points in $J$ .

\subsection{The toric variety $B_{3}$}

\ This variety  is defined by a fan $\Delta$ such that elements in
$\Delta(1)$ satisfy :
\begin{equation}
\left\{ \begin{aligned}
         \ v_{1}+v_{2}+v_{3}+v_{4}=v_{6} \\
                  \ v_{5}+v_{6}=0
                          \end{aligned} \right. \notag
\end{equation}

$\mathbb{C}^{*}\times \mathbb{C}^{*}$ acts on this toric variety as
follows :

\begin{equation}
(X_{1}, \ldots , X_{6} )\rightarrow
(\lambda_{1}X_{1},\lambda_{1}X_{2}, \lambda_{1}X_{3},
\lambda_{1}X_{4},\lambda_{1}\lambda_{2}X_{5},
\lambda_{2}X_{6})\notag
\end{equation}

Under this action , $X\simeq \mathbb{C}^{6}-Z\diagup
\mathbb{C}^{*}\times \mathbb{C}^{*}$ , where
$Z=\{X_{1}X_{5}=X_{2}X_{5}=X_{3}X_{5}=X_{4}X_{5}
=X_{1}X_{6}=X_{2}X_{6}=X_{3}X_{6}=X_{4}X_{6}=0\}$

The anti-canonical forms of $X$ are linear combinations of the
following forms :

\begin{equation}
\begin{split}
 X_{5}^{2}f_{3}(X_{1},\ldots ,
X_{4}), X_{5}X_{6}f_{4}(X_{1},\ldots , X_{4}) ,
\\
X_{6}^{2}f_{5}(X_{1},\ldots , X_{4}).
\end{split}  \notag
\end{equation}
where $f_{i}(X_{1},\ldots , X_{4})$ are degree $i$ homogenous forms
, for $i\geq 1$ .

 Now $H^{2}(X, \mathbb{Z})\simeq H_{6}(X ,
\mathbb{Z})\simeq A_{3}(X)$ is a rank $2$ free abelian group , and
$[D_{4}] , [D_{5}]$ form a base of this group .

\begin{table}[!h]
\tabcolsep 0.5mm \caption{$B_{3}$}
\begin{center}
\begin{tabular}{r@{}lr@{}lr@{}lr@{}l}
\hline \multicolumn{2}{c}{rational curve
$C_{i0}$}&\multicolumn{2}{c}{type}&\multicolumn{2}{c}{cohomology
class}&\multicolumn{2}{c}{anti-canonical section $Y_{i0}$}
\\ \hline
(s,t,0,0,1,0&)   &(1,1,1,1,0,-1&)&[$D_{4}$&]&$X_{3}X_{1}^{2}X_{5}^{2}+X_{4}X_{2}^{2}X_{5}^{2}+X_{6}X_{5}X_{2}^{4}=$&0 \\
(s,t,0,0,0,1&)&(1,1,1,1,1,0&)&[$D_{4}]+[D_{5}$&]&$X_{3}X_{1}^{4}X_{6}^{2}+X_{4}X_{1}^{2}X_{2}^{2}X_{6}^{2}+X_{5}X_{6}X_{2}^{4}=$&0\\
(1,0,0,0,s,t&)&(0,0,0,0,1,1&) &[$D_{5}$&]  & $X_{2}X_{1}^{2}X_{5}^{2}+X_{3}X_{1}^{3}X_{5}X_{6}+X_{4}X_{1}^{4}X_{6}^{2}=$&0\\
\hline
\end{tabular}
\end{center}
\end{table}

Denote the three rational curves in the first column by $C_{10}$ ,
$C_{20}$ , $C_{30}$ respectively , and denote the corresponding
parameterizing spaces by $\mathcal {M}_{1}$ , $\mathcal {M}_{2}$ ,
$\mathcal {M}_{3}$ , then take the corresponding three Zariski open
sets $U_{1}$ , $U_{2}$ , $U_{3}$ such that the rational curves
parameterized by points in $U_{i}$ are smooth  , and  the homogenous
coordinates $X_{5}$ , $X_{6}$ , $X_{1}$ have no zero points on
rational curves in $U_{1}$ , $U_{2}$ , $U_{3}$ respectively .

The case of $C_{1}$ and $C_{2}$

This case is trivial . Since the homogenous coordinate $X_{6}$ has
no zero points on $C_{2}$ , and it's always zero on $C_{1}$ , so
$U_{12}=\emptyset$.

The case of $C_{1}$ and $C_{3}$

$dim$ $U_{13}= 11$ , $C_{1}$ , $C_{3}$ can be assumed to be $
(s,t,0,0,1,0)$ , $( 1,0, 0,0,s,t )$ , and the following $6$ forms
will suffice :
\begin{equation}
X_{1}^{4}X_{5}X_{6},
X_{1}^{5}X_{6}^{2},X_{5}^{2}X_{1}^{i}X_{2}^{3-i}(0\leq i\leq 3)
.\notag
\end{equation}

The case of $C_{2}$ and $C_{3}$

$dim$ $U_{23}=13$ , $C_{2}$ , $C_{3}$ can be assumed to be
$(s,t,0,0,\alpha_{1}(s,t),1)$ , $(1,0,0,0, s,t)$ , and the following
$8$ forms will suffice :
\begin{equation}
X_{1}^{4}X_{5}X_{6},
X_{1}^{5}X_{6}^{2},X_{5}^{2}X_{1}^{3},X_{6}^{2}X_{2}X_{1}^{i}X_{2}^{4-i}(0\leq
i\leq 4) .\notag
\end{equation}

\subsection{The toric variety $B_{4}$}

This variety  is defined by a fan $\Delta$ such that elements in
$\Delta(1) $ satisfy :
 \begin{equation}
\left\{ \begin{aligned}
         \ v_{1}+v_{2}+v_{3}+v_{4}=0 \\
                  \ v_{5}+v_{6}=0
                          \end{aligned} \right. \notag
\end{equation}

$\mathbb{C}^{*}\times \mathbb{C}^{*}$ acts on this toric variety as
follows :

\begin{equation}
(X_{1}, \ldots , X_{6} )\rightarrow
(\lambda_{1}X_{1},\lambda_{1}X_{2}, \lambda_{1}X_{3},
\lambda_{1}X_{4},\lambda_{2}X_{5}, \lambda_{2}X_{6})\notag
\end{equation}

Under this action , $X\simeq \mathbb{C}^{6}-Z\diagup
\mathbb{C}^{*}\times \mathbb{C}^{*}\simeq
\mathbb{P}^{1}\times\mathbb{P}^{3}$ , where
$Z=\{X_{1}X_{5}=X_{2}X_{5}=X_{3}X_{5}=X_{4}X_{5}
=X_{1}X_{6}=X_{2}X_{6}=X_{3}X_{6}=X_{4}X_{6}=0\}$

The   anti-canonical forms of $X$ are linear combinations of  the
forms $f_{4}(X_{1},X_{2},X_{3},X_{4})g_{2}(X_{5},X_{6})$ ,  where
$f_{i}(X_{1},X_{2} , X_{3},X_{4})$ , $g_{i}(X_{5},X_{6})$ are degree
$i$ homogenous forms  , for $i\geq 1$ .

 Now $H^{2}(X, \mathbb{Z})\simeq H_{6}(X ,
\mathbb{Z})\simeq A_{3}(X)$ is a rank $2$ free abelian group , and
$[D_{4}]$ , $ [D_{5}]$ form a base of this group .

\begin{table}[!h]
\tabcolsep 0.5mm \caption{$B_{4}$}
\begin{center}
\begin{tabular}{r@{}lr@{}lr@{}lr@{}l}
\hline \multicolumn{2}{c}{rational curve
$C_{i0}$}&\multicolumn{2}{c}{type}&\multicolumn{2}{c}{cohomology
class}&\multicolumn{2}{c}{anti-canonical section $Y_{i0}$}
\\ \hline
(1,0,0,0,s,t&)   &(0,0,0,0,1,1&)&[$D_{5}$&]&$X_{2}X_{5}^{2}X_{1}^{3}+X_{3}X_{5}X_{6}X_{1}^{3}+X_{4}X_{6}^{2}X_{1}^{3}=$&0 \\
(s,t,0,0,1,0&)&(1,1,1,1,0,0&)&[$D_{4}$&]&$X_{3}X_{5}^{2}X_{1}^{3}+X_{4}X_{5}^{2}X_{2}^{3}+X_{6}X_{1}^{2}X_{5}X_{2}^{2}=$&0\\
(s,t,0,0,s,t&)&(1,1,1,1,1,1&) &[$D_{4}$]+[$D_{5}$&]  & $X_{3}X_{1}^{3}X_{5}^{2}+X_{4}X_{1}^{3}X_{6}^{2}+X_{2}^{4}X_{5}X_{6}-X_{1}X_{2}^{3}X_{6}^{2}=$&0\\
\hline
\end{tabular}
\end{center}
\end{table}

Denote the three rational curves in the first column by $C_{10}$ ,
$C_{20}$ , $C_{30}$ respectively , and denote the corresponding
parameterizing spaces by $\mathcal {M}_{1}$ , $\mathcal {M}_{2}$ ,
$\mathcal {M}_{3}$ , then take the corresponding three Zariski open
sets $U_{1}$ , $U_{2}$ , $U_{3}$ such that the rational curves
parameterized by points in $U_{i}$ are smooth  .

The case of $C_{1}$ and $C_{2}$

$dim$ $U_{12}=12$ , $C_{1}$ , $C_{2}$ can be assumed to be
$(1,0,0,0,s,t)$ , $(s,t,0,0, 1,0)$ , and the following $7$ forms
will suffice :
\begin{equation}
X_{5}^{i}X_{6}^{4-i}X_{1}^{2}(0\leq i\leq 4), X_{1}X_{2}X_{5}^{4},
X_{2}^{2}X_{5}^{4} . \notag
\end{equation}

The case of $C_{1}$ and $C_{3}$

$dim$ $U_{13}= 14$ , $C_{1}$ , $C_{3}$ can be assumed to be $
(1,0,0,0,s,t)$ , $( s,t,0,s,0,t)$ or \\  $(1,0,0,0,s,t)$ , $(
s,t,0,0,s,t)$ . In the first case ,  the following $10$ forms will
suffice :
\begin{equation}
X_{5}^{i}X_{6}^{4-i}X_{1}^{2}(0\leq i\leq
4),X_{4}^{i}X_{6}^{4-i}X_{2}^{2}(0\leq i\leq 4) .\notag
\end{equation}

And in the second case , the following $10$ forms will suffice :

\begin{equation}
X_{5}^{i}X_{6}^{4-i}X_{1}^{2}(0\leq i\leq
4),X_{5}^{i}X_{6}^{4-i}X_{2}^{2}(0\leq i\leq 4) .\notag
\end{equation}

The case of $C_{2}$ and $C_{3}$

$dim$ $U_{12}=11$ , $C_{2}$ , $C_{3}$ can be assumed to be
$(s,t,0,0,1,0)$ , $(s,t,0,0, s,t)$  , and the following $6$ forms
will suffice :
\begin{equation}
X_{1}^{i}X_{2}^{2-i}X_{5}^{4}(0\leq i\leq
2),X_{1}^{i}X_{2}^{2-i}X_{6}^{4}(0\leq i\leq 2).\notag
\end{equation}

\subsection{The toric variety $B_{5}$}

This variety  is defined by a fan $\Delta$ such that elements in
$\Delta(1) $ satisfy :
\begin{equation}
\left\{ \begin{aligned}
         \ v_{1}+v_{2}+v_{3}+v_{4}=0 \\
                  \ v_{5}+v_{6}=v_{4}
                          \end{aligned} \right.\notag
\end{equation}

$\mathbb{C}^{*}\times \mathbb{C}^{*}$ acts on this toric variety as
follows :

\begin{equation}
(X_{1}, \ldots , X_{6} )\rightarrow
(\lambda_{1}\lambda_{2}X_{1},\lambda_{1}\lambda_{2}X_{2},
\lambda_{1}\lambda_{2}X_{3}, \lambda_{1}X_{4},\lambda_{2}X_{5},
\lambda_{2}X_{6})\notag
\end{equation}

Under this action , $X\simeq \mathbb{C}^{6}-Z\diagup
\mathbb{C}^{*}\times \mathbb{C}^{*}$ , where
$Z=\{X_{1}X_{5}=X_{2}X_{5}=X_{3}X_{5}=X_{4}X_{5}
=X_{1}X_{6}=X_{2}X_{6}=X_{3}X_{6}=X_{4}X_{6}=0\}$

The anti-canonical forms of $X$ are linear combinations of the
following forms :

\begin{equation}
\begin{split}
X_{4}^{4}f_{5}(X_{5},X_{6}) , X_{4}^{3}g_{1}(X_{1},X_{2} ,
X_{3})f_{4}(X_{5}, X_{6}), X_{4}^{2}g_{2}(X_{1},X_{2} ,
X_{3})f_{3}(X_{5}, X_{6}),\\
X_{4}g_{3}(X_{1},X_{2} , X_{3})f_{2}(X_{5}, X_{6}),
g_{4}(X_{1},X_{2}, X_{3})f_{1}(X_{5}, X_{6})
\end{split}  \notag
\end{equation}
where $g_{i}(X_{1},X_{2} , X_{3})$ , $f_{i}(X_{5},X_{6})$ are degree
$i$ homogenous forms  , for $i\geq 1$ .

 Now $H^{2}(X, \mathbb{Z})\simeq H_{6}(X ,
\mathbb{Z})\simeq A_{3}(X)$ is a rank $2$ free abelian group , and
$[D_{4}]$ , $[D_{5}]$ form a base of this group .

\begin{table}[!h]
\tabcolsep 0.5mm \caption{$B_{5}$}
\begin{center}
\begin{tabular}{r@{}lr@{}lr@{}lr@{}l}
\hline \multicolumn{2}{c}{rational curve
$C_{i0}$}&\multicolumn{2}{c}{type}&\multicolumn{2}{c}{cohomology
class}&\multicolumn{2}{c}{anti-canonical section $Y_{i0}$}
\\ \hline
(0,0,0,1,s,t&)   &(1,1,1,0,1,1&)&[$D_{5}$&]&$X_{1}X_{4}^{3}X_{6}^{4}+X_{2}X_{4}^{3}X_{5}^{4}+X_{3}X_{4}^{3}X_{5}^{2}X_{6}^{2}=$&0 \\
(s,0,0,t,1,0&)&(1,1,1,1,0,0&)&[$D_{4}$&]&$X_{2}X_{1}^{3}X_{5}+
X_{3}X_{4}^{2}X_{1}X_{5}^{3}+X_{6}X_{5}^{4}X_{4}^{4}=$&0\\
(1,0,0,0,s,t&)&(0,0,0,-1,1,1&) &-[$D_{4}$]+[$D_{5}$&]  & $X_{2}X_{1}^{3}X_{5}+X_{3}X_{1}^{3}X_{6}+X_{4}X_{6}^{2}X_{1}^{3}=$&0\\
\hline
\end{tabular}
\end{center}
\end{table}

Denote the three rational curves in the first column by $C_{10}$ ,
$C_{20}$ , $C_{30}$ respectively , and denote the corresponding
parameterizing spaces by $\mathcal {M}_{1}$ , $\mathcal {M}_{2}$ ,
$\mathcal {M}_{3}$ , then take the corresponding three Zariski open
sets $U_{1}$ , $U_{2}$ , $U_{3}$ such that the rational curves
parameterized by points in $U_{i}$ are smooth , and  the homogenous
coordinates $X_{4}$ , $X_{5}$ , $X_{1}$ have no zero points on
rational curves in $U_{1}$ , $U_{2}$ , $U_{3}$ respectively .


The case of $C_{1}$ and $C_{2}$

$dim$ $U_{12}= 15$ , $C_{1}$ , $C_{2}$ can be assumed to be
$(\alpha_{1}(s,t),\beta_{1}(s,t),\gamma_{1}(s,t),1,s,t)$ , $(s,
\tilde{\alpha}_{1}(s,t), \tilde{\beta}_{1}(s,t), t, 1, 0)$ , and the
following $10 $ forms will suffice :
\begin{equation}
X_{4}^{4}X_{5}^{5}, X_{4}^{3}X_{1}X_{5}^{4},
X_{4}^{2}X_{1}^{2}X_{5}^{3}, X_{4}X_{1}^{3}X_{5}^{2},
X_{1}^{4}X_{5}, X_{4}^{4}X_{6}X_{5}^{i}X_{6}^{4-i}(0\leq i\leq
4).\notag
\end{equation}

The case of $C_{2}$ and $C_{3}$

$dim$ $U_{23}= 11$ , $C_{2}$ , $C_{3}$ can be assumed to be $ (s,
\tilde{\alpha}_{1}(s,t),\tilde{\beta}_{1}(s,t), t, 1, 0)$ ,
$(1,0,0,0,s,t)$ , and the following $6$ forms will suffice :
\begin{equation}
X_{4}^{4}X_{5}^{5}, X_{4}^{3}X_{1}X_{5}^{4},
X_{4}^{2}X_{1}^{2}X_{5}^{3}, X_{4}X_{1}^{3}X_{5}^{2},
X_{1}^{4}X_{5}, X_{1}^{4}X_{6} .\notag
\end{equation}

The case of $C_{1}$ and $C_{3}$

This case is trivial . Since the homogenous coordinate $X_{4}$ has
no zero points on $C_{1}$ , and it's always zero on $C_{3}$ , so
$U_{13}=\emptyset$.

\subsection{The toric variety $C_{1}$}

This variety is defined by a fan $\Delta$ such that elements in
$\Delta(1)$ satisfy :
\begin{equation}
\left\{ \begin{aligned}
         \ v_{1}+v_{2}+v_{3}=0 \\
                  \ v_{4}+v_{5}+v_{6}=2v_{3}
                          \end{aligned} \right.\notag
\end{equation}

$\mathbb{C}^{*}\times \mathbb{C}^{*}$ acts on this toric variety as
follows :

\begin{equation}
(X_{1}, \ldots , X_{6} )\rightarrow
(\lambda_{1}\lambda_{2}^{2}X_{1},\lambda_{1}\lambda_{2}^{2}X_{2},
\lambda_{1}X_{3}, \lambda_{2}X_{4},\lambda_{2}X_{5},
\lambda_{2}X_{6})\notag
\end{equation}

Under this action , $X\simeq \mathbb{C}^{6}-Z\diagup
\mathbb{C}^{*}\times \mathbb{C}^{*}$ , where
$Z=\{X_{1}X_{4}=X_{1}X_{5}=X_{1}X_{6}=X_{2}X_{4}
=X_{2}X_{5}=X_{2}X_{6}=X_{3}X_{4}=X_{3}X_{5}=X_{3}X_{6}=0\}$

The anti-canonical forms of $X$ are linear combinations of the
following forms :

\begin{gather}
X_{3}^{3}f_{7}(X_{4},X_{5},X_{6}) ,
X_{3}^{2}f_{5}(X_{4},X_{5},X_{6})g_{1}(X_{1}, X_{2}),
X_{3}f_{3}(X_{4},X_{5}, X_{6})g_{2}(X_{1} ,
X_{2}),\notag \\
f_{1}(X_{4},X_{5}, X_{6})g_{3}(X_{1}, X_{2}).\notag
\end{gather}
where $g_{i}(X_{1}, X_{2})$ , $f_{i}(X_{4}, X_{5},X_{6})$ are degree
$i$ homogenous forms  , for $i\geq 1$ .

 Now $H^{2}(X, \mathbb{Z})\simeq H_{6}(X ,
\mathbb{Z})\simeq A_{3}(X)$ is a rank $2$ free abelian group , and
$[D_{3}]$ , $[D_{6}]$ form a base of this group .

\begin{table}[!h]
\tabcolsep 0.5mm \caption{$C_{1}$}
\begin{center}
\begin{tabular}{r@{}lr@{}lr@{}lr@{}l}
\hline \multicolumn{2}{c}{rational curve
$C_{i0}$}&\multicolumn{2}{c}{type}&\multicolumn{2}{c}{cohomology
class}&\multicolumn{2}{c}{anti-canonical section $Y_{i0}$}
\\ \hline
(s,0,t,0,0,1&)   &(1,1,1,0,0,0&)&[$D_{3}$&]&$X_{2}X_{1}^{2}X_{6}+X_{4}X_{3}^{2}X_{6}^{4}X_{1}+X_{5}X_{3}^{3}X_{6}^{6}=$&0 \\
(0,0,1,s,0,t&)&(2,2,0,1,1,1&)&[$D_{6}$&]&$X_{1}X_{3}^{2}X_{4}^{5}+
X_{2}X_{4}^{2}X_{3}^{2}X_{6}^{3}+X_{5}X_{3}^{3}X_{6}^{6}=$&0\\
(s,t,0,s,0,t&)&(1,1,-1,1,1,1&) &[$D_{6}$]-[$D_{3}$&]  & $X_{3}X_{2}^{2}X_{6}^{3}+X_{5}X_{1}^{3}+X_{1}X_{6}X_{2}^{2}-X_{2}^{3}X_{4}=$&0\\
\hline
\end{tabular}
\end{center}
\end{table}

Denote the three rational curves in the first column by $C_{10}$ ,
$C_{20}$ , $C_{30}$ respectively , and denote the corresponding
parameterizing spaces by $\mathcal {M}_{1}$ , $\mathcal {M}_{2}$ ,
$\mathcal {M}_{3}$ , then take the corresponding three Zariski open
sets $U_{1}$ , $U_{2}$ , $U_{3}$ such that the rational curves
parameterized by points in $U_{i}$ are smooth , and the homogenous
coordinates $X_{6}$ , $ X_{3}$ have no zero points on rational
curves in $U_{1}$ , $U_{2}$ respectively .

The case of $C_{1}$ and $C_{2}$

$dim$ $U_{12}=16$ , $C_{1}$ , $C_{2}$ can be assumed to be
$(s,\alpha_{1}(s,t),t,0,0,1)$ ,
$(\alpha_{2}(s,t),\beta_{2}(s,t),1,s,0, t)$ , and the following $11
$ forms will suffice :
\begin{equation}
X_{3}^{3}X_{6}^{7}, X_{3}^{2}X_{1}X_{6}^{5},
X_{1}^{2}X_{6}^{3}X_{3}, X_{6}X_{1}^{3},
X_{3}^{3}X_{4}X_{4}^{i}X_{6}^{6-i}(0\leq i\leq 6) .\notag
\end{equation}

The case of $C_{1}$ and $C_{3}$

$dim$ $U_{13}= 13$ , $C_{1}$ , $C_{3}$ can be assumed to be $
(s,\alpha_{1}(s,t),t,0,0,1)$ , $(
\tilde{\alpha}_{1}(s,t),\tilde{\beta}_{1}(s,t), 0,s,0,t )$ , and the
following $8$ forms will suffice :
\begin{equation}
X_{3}^{3}X_{6}^{7}, X_{3}^{2}X_{1}X_{6}^{5},
X_{1}^{2}X_{6}^{3}X_{3}, X_{6}X_{1}^{3}, X_{4}X_{1}^{3},
X_{4}X_{1}^{2}X_{2}, X_{4}X_{1}X_{2}^{2}, X_{4}X_{2}^{3} .\notag
\end{equation}

The case of $C_{2}$ and $C_{3}$

This case is trivial . Since the homogenous coordinate $X_{3}$ has
no zero points on $C_{2}$ , and it's always zero on $C_{3}$ , so
$U_{23}=\emptyset$ .

\subsection{The toric variety $C_{2}$}

This variety is defined by a fan $\Delta$ such that elements in
$\Delta(1)$ satisfy :
\begin{equation}
\left\{ \begin{aligned}
         \ v_{1}+v_{2}+v_{3}=0 \\
                  \ v_{4}+v_{5}+v_{6}=v_{3}
                          \end{aligned} \right.\notag
\end{equation}

$\mathbb{C}^{*}\times \mathbb{C}^{*}$ acts on this toric variety as
follows :

\begin{equation}
(X_{1}, \ldots , X_{6} )\rightarrow
(\lambda_{1}\lambda_{2}X_{1},\lambda_{1}\lambda_{2}X_{2},
\lambda_{1}X_{3}, \lambda_{2}X_{4},\lambda_{2}X_{5},
\lambda_{2}X_{6})\notag
\end{equation}

Under this action , $X\simeq \mathbb{C}^{6}-Z\diagup
\mathbb{C}^{*}\times \mathbb{C}^{*}$ , where
$Z=\{X_{1}X_{4}=X_{1}X_{5}=X_{1}X_{6}=X_{2}X_{4}
=X_{2}X_{5}=X_{2}X_{6}=X_{3}X_{4}=X_{3}X_{5}=X_{3}X_{6}=0\}$

The anti-canonical forms of $X$ are linear combinations of the
following forms :

\begin{gather}
X_{3}^{3}f_{5}(X_{4},X_{5},X_{6}) ,
X_{3}^{2}f_{4}(X_{4},X_{5},X_{6})g_{1}(X_{1}, X_{2}),
X_{3}f_{3}(X_{4},X_{5}, X_{6})g_{2}(X_{1} ,
X_{2}),\notag \\
f_{2}(X_{4},X_{5}, X_{6})g_{3}(X_{1}, X_{2}) . \notag
\end{gather}
where $g_{i}(X_{1}, X_{2})$ , $f_{i}(X_{4},X_{5},X_{6})$ are degree
$i$ homogenous forms  , for $i\geq 1$ .

 Now $H^{2}(X, \mathbb{Z})\simeq H_{6}(X ,
\mathbb{Z})\simeq A_{3}(X)$ is a rank $2$ free abelian group , and
$[D_{3}]$ , $[D_{6}]$ form  a base of this group .

\begin{table}[!h]
\tabcolsep 0.5mm \caption{$C_{2}$}
\begin{center}
\begin{tabular}{r@{}lr@{}lr@{}lr@{}l}
\hline \multicolumn{2}{c}{rational curve
$C_{i0}$}&\multicolumn{2}{c}{type}&\multicolumn{2}{c}{cohomology
class}&\multicolumn{2}{c}{anti-canonical section $Y_{i0}$}
\\ \hline
(s,0,t,0,0,1&)   &(1,1,1,0,0,0&)&[$D_{3}$&]&$X_{2}X_{1}^{2}X_{6}^{2}+X_{4}X_{3}^{2}X_{6}^{3}X_{1}+X_{5}X_{3}^{3}X_{6}^{4}=$&0 \\
(0,0,1,s,0,t&)&(1,1,0,1,1,1&)&[$D_{6}$&]&$X_{1}X_{3}^{2}X_{4}^{4}+
X_{2}X_{4}^{2}X_{3}^{2}X_{6}^{2}+X_{5}X_{3}^{3}X_{6}^{4}=$&0\\
(1,0,0,s,0,t&)&(0,0,-1,1,1,1&) &-[$D_{3}$]+[$D_{6}$&]  & $X_{2}X_{1}^{2}X_{4}^{2}+X_{3}X_{2}^{2}X_{5}^{3}+X_{5}X_{1}^{3}X_{6}=$&0\\
\hline
\end{tabular}
\end{center}
\end{table}

Denote the three rational curves in the first column by $C_{10}$ ,
$C_{20}$ , $C_{30}$ respectively , and denote the corresponding
parameterizing spaces by $\mathcal {M}_{1}$ , $\mathcal {M}_{2}$ ,
$\mathcal {M}_{3}$ , then take the corresponding three Zariski open
sets $U_{1}$ , $U_{2}$ , $U_{3}$ such that the rational curves
parameterized by points in $U_{i}$ are smooth , and  the homogenous
coordinates $X_{6}$ , $ X_{3}$ have no zero points on rational
curves in $U_{1}$ , $U_{2}$ respectively .

The case of $C_{1}$ and $C_{2}$

$dim$ $U_{12}=14$ , $C_{1}$ , $C_{2}$ can be assumed to be
$(s,\alpha_{1}(s,t),t,0,0,1)$ , $
(\tilde{\alpha}_{1}(s,t),\tilde{\beta}_{1}(s,t),1,s,0, t)$ , and the
following $9 $ forms will suffice :
\begin{equation}
X_{3}^{3}X_{6}^{5}, X_{3}^{2}X_{1}X_{6}^{4},
X_{1}^{2}X_{6}^{3}X_{3}, X_{6}^{2}X_{1}^{3},
X_{3}^{3}X_{4}X_{4}^{i}X_{6}^{4-i}(0\leq i\leq 4).\notag
\end{equation}

The case of $C_{1}$ and $C_{3}$

$dim$ $U_{13}= 11$ , $C_{1}$ , $C_{3}$ can be assumed to be $
(s,\alpha_{1}(s,t),t,0,0,1)$ , $( 1,0, 0,s,0,t )$ , and the
following $6$ forms will suffice :
\begin{equation}
X_{3}^{3}X_{6}^{5}, X_{3}^{2}X_{1}X_{6}^{4},
X_{1}^{2}X_{6}^{3}X_{3}, X_{6}^{2}X_{1}^{3},X_{4}X_{6}X_{1}^{3},
X_{4}^{2}X_{1}^{3}.\notag
\end{equation}

The case of $C_{2}$ and $C_{3}$

This case is trivial . Since the homogenous coordinate $X_{3}$ has
no zero points on $C_{2}$ , and it's always zero on $C_{3}$ , so
$U_{23}=\emptyset$ .

\subsection{The toric variety $C_{3}$}

This variety  is defined by a fan $\Delta$ such that elements in
$\Delta(1)$ satisfy :
\begin{equation}
\left\{ \begin{aligned}
         \ v_{1}+v_{2}+v_{3}=0 \\
                  \ v_{4}+v_{5}+v_{6}=v_{1}+v_{2}=-v_{3}
                          \end{aligned} \right.\notag
\end{equation}

$\mathbb{C}^{*}\times \mathbb{C}^{*}$ acts on this toric variety as
follows :

\begin{equation}
(X_{1}, \ldots , X_{6} )\rightarrow
(\lambda_{1}\lambda_{2}^{-1}X_{1},\lambda_{1}\lambda_{2}^{-1}X_{2},
\lambda_{1}X_{3}, \lambda_{2}X_{4},\lambda_{2}X_{5},
\lambda_{2}X_{6})\notag
\end{equation}

Under this action , $X\simeq \mathbb{C}^{6}-Z\diagup
\mathbb{C}^{*}\times \mathbb{C}^{*}$ , where
$Z=\{X_{1}X_{4}=X_{1}X_{5}=X_{1}X_{6}=X_{2}X_{4}
=X_{2}X_{5}=X_{2}X_{6}=X_{3}X_{4}=X_{3}X_{5}=X_{3}X_{6}=0\}$

The anti-canonical forms of $X$ are linear combinations of the
following forms :

\begin{gather}
X_{3}^{3}f_{1}(X_{4},X_{5},X_{6}) ,
X_{3}^{2}f_{2}(X_{4},X_{5},X_{6})g_{1}(X_{1}, X_{2}),
X_{3}f_{3}(X_{4},X_{5}, X_{6})g_{2}(X_{1} ,
X_{2}),\notag \\
f_{4}(X_{4},X_{5}, X_{6})g_{3}(X_{1}, X_{2}).\notag
\end{gather}
where $g_{i}(X_{1},X_{2})$ , $f_{i}(X_{4},X_{5},X_{6})$ are degree
$i$ homogenous forms  , for $i\geq 1$ .

 Now $H^{2}(X, \mathbb{Z})\simeq H_{6}(X ,
\mathbb{Z})\simeq A_{3}(X)$ is a rank $2$ free abelian group , and
$[D_{3}]$ , $[D_{6}]$ form a base of this group .

\begin{table}[!h]
\tabcolsep 0.5mm \caption{$C_{3}$}
\begin{center}
\begin{tabular}{r@{}lr@{}lr@{}lr@{}l}
\hline \multicolumn{2}{c}{rational curve
$C_{i0}$}&\multicolumn{2}{c}{type}&\multicolumn{2}{c}{cohomology
class}&\multicolumn{2}{c}{anti-canonical section $Y_{i0}$}
\\ \hline
(s,0,t,0,0,1&)   &(1,1,1,0,0,0&)&[$D_{3}$&]&$X_{2}X_{1}^{2}X_{6}^{4}+X_{4}X_{3}^{2}X_{6}X_{1}+X_{5}X_{3}^{3}=$&0 \\
(0,0,1,s,0,t&)&(-1,-1,0,1,1,1&)&[$D_{6}$&]&$X_{1}X_{3}^{2}X_{4}^{2}+
X_{2}X_{3}^{2}X_{4}^{2}+X_{5}X_{3}^{3}=$&0\\
(1,0,0,s,0,t&)&(0,0,1,1,1,1&) &[$D_{3}$]+[$D_{6}$&]  & $X_{2}X_{1}^{2}X_{4}^{4}+X_{3}X_{1}^{2}X_{4}^{2}X_{6}+X_{5}X_{1}^{3}X_{6}^{3}=$&0\\
\hline
\end{tabular}
\end{center}
\end{table}

Denote the three rational curves in the first column by $C_{10}$ ,
$C_{20}$ , $C_{30}$ respectively , and denote the corresponding
parameterizing spaces by $\mathcal {M}_{1}$ , $\mathcal {M}_{2}$ ,
$\mathcal {M}_{3}$ , then take the corresponding three Zariski open
sets $U_{1}$ , $U_{2}$ , $U_{3}$ such that the rational curves
parameterized by points in $U_{i}$ are smooth , and  the homogenous
coordinates $X_{6}$ , $X_{3}$ , $X_{1}$ have no zero points on
rational curves in $U_{1}$ , $U_{2}$ , $U_{3}$ respectively .

The case of $C_{1}$ and $C_{2}$

$dim$ $U_{12}=10$ , $C_{1}$ , $C_{2}$ can be assumed to be
$(s,\alpha_{1}(s,t),t,0,0,1)$ , $ (0,0,1,s,0, t)$ , and the
following $5$ forms will suffice :
\begin{equation}
X_{3}^{3}X_{6}, X_{3}^{2}X_{1}X_{6}^{2}, X_{1}^{2}X_{6}^{3}X_{3},
X_{6}^{4}X_{1}^{3}, X_{3}^{3}X_{4}.\notag
\end{equation}

The case of $C_{1}$ and $C_{3}$

$dim$ $U_{13}= 13$ , $C_{1}$ , $C_{3}$ can be assumed to be $
(s,\alpha_{1}(s,t),t,0,0,1)$ , $( 1,0, \tilde{\alpha}_{1}(s,t),t,0,s
)$ , and the following $8$ forms will suffice :
\begin{equation}
X_{3}^{3}X_{6}, X_{3}^{2}X_{1}X_{6}^{2}, X_{1}^{2}X_{6}^{3}X_{3},
X_{6}^{4}X_{1}^{3}, X_{4}X_{1}^{3}X_{4}^{i}X_{6}^{3-i}(0\leq i\leq
3).\notag
\end{equation}

The case of $C_{2}$ and $C_{3}$

This case is trivial . Since the homogenous coordinate $X_{1}$ has
no zero points on $C_{3}$ , and it's always zero on $C_{2}$ , so
$U_{23}=\emptyset$ .

\subsection{The toric variety $C_{4}$}

This variety is defined by a fan $\Delta$ such that elements in
$\Delta(1)$ satisfy :
\begin{equation}
\left\{ \begin{aligned}
         \ v_{1}+v_{2}+v_{3}=0 \\
                  \ v_{4}+v_{5}+v_{6}=0
                          \end{aligned} \right. \notag
\end{equation}

$\mathbb{C}^{*}\times \mathbb{C}^{*}$ acts on this toric variety as
follows :

\begin{equation}
(X_{1}, \ldots , X_{6} )\rightarrow
(\lambda_{1}X_{1},\lambda_{1}X_{2}, \lambda_{1}X_{3},
\lambda_{2}X_{4},\lambda_{2}X_{5}, \lambda_{2}X_{6})\notag
\end{equation}

Under this action , $X\simeq \mathbb{C}^{6}-Z\diagup
\mathbb{C}^{*}\times \mathbb{C}^{*}\simeq \mathbb{P}^{2}\times
\mathbb{P}^{2}$ , where
$Z=\{X_{1}X_{4}=X_{1}X_{5}=X_{1}X_{6}=X_{2}X_{4}
=X_{2}X_{5}=X_{2}X_{6}=X_{3}X_{4}=X_{3}X_{5}=X_{3}X_{6}=0\}$

The anti-canonical forms of $X$ are linear combinations of  the
forms $f_{3}(X_{1},X_{2},X_{3})g_{3}(X_{4},X_{5},X_{6})$ , where
$f_{i}(X_{1},X_{2} , X_{3})$ , $g_{i}(X_{4},X_{5},X_{6})$ are degree
$i$ homogenous forms  , for $i\geq 1$ .

 Now $H^{2}(X, \mathbb{Z})\simeq H_{6}(X ,
\mathbb{Z})\simeq A_{3}(X)$ is a rank $2$ free abelian group , and
$[D_{3}]$ , $[D_{6}]$ form a base of this group .

\begin{table}[!h]
\tabcolsep 0.2mm \caption{$C_{4}$}
\begin{center}
\begin{tabular}{r@{}lr@{}lr@{}lr@{}l}
\hline \multicolumn{2}{c}{rational curve
$C_{i0}$}&\multicolumn{2}{c}{type}&\multicolumn{2}{c}{cohomology
class}&\multicolumn{2}{c}{anti-canonical section $Y_{i0}$}
\\ \hline
(s,t,0,1,0,0&)   &(1,1,1,0,0,0&)&[$D_{3}$&]&$X_{5}X_{4}^{2}X_{2}^{3}+X_{6}X_{4}^{2}X_{1}X_{2}^{2}+X_{3}X_{4}^{3}X_{1}^{2}=$&0 \\
(1,0,0,s,t,0&)&(0,0,0,1,1,1&)&[$D_{6}$&]&$X_{2}X_{1}^{2}X_{5}^{3}+X_{3}X_{1}^{2}X_{4}X_{5}^{2}+X_{6}X_{1}^{3}X_{4}^{2}=$&0\\
(s,t,0,s,t,0&)&(1,1,1,1,1,1&) & [$D_{3}$]+[$D_{6}$&]  & $X_{3}X_{1}^{2}X_{4}X_{5}^{2}+X_{6}X_{1}^{3}X_{4}^{2}+X_{1}X_{2}^{2}X_{5}^{3}-X_{2}^{3}X_{4}X_{5}^{2}=$&0\\
\hline
\end{tabular}
\end{center}
\end{table}

In the following Table $9$ , denote the three rational curves in the
first column by $C_{10}$ , $C_{20}$ , $C_{30}$ respectively , and
denote the corresponding parameterizing spaces by $\mathcal {M}_{1}$
, $\mathcal {M}_{2}$ , $\mathcal {M}_{3}$ , then take the
corresponding three Zariski open sets $U_{1}$ , $U_{2}$ , $U_{3}$
such that the rational curves parameterized by points in $U_{i}$ are
smooth .

The case of $C_{1}$ and $C_{2}$

$dim$ $U_{12}=12$ , $C_{1}$ , $C_{2}$ can be assumed to be
$(s,t,0,1,0,0)$ , $ (1,0,0,s,t, 0)$ , and the following $7$ forms
will suffice :
\begin{equation}
X_{1}^{i}X_{2}^{3-i}X_{4}^{3}(0\leq i\leq 3),
X_{1}^{3}X_{5}X_{4}^{i}X_{5}^{2-i}(0\leq i\leq 2).\notag
\end{equation}

The case of $C_{1}$ and $C_{3}$

$dim$ $U_{13}= 13$ . $C_{1}$ , $C_{3}$ can be assumed to be $
(s,t,0,1,0,0)$ , $( s,0,t,s,t,0)$ or  $(s,t,0,1,0,0)$ , $(
s,t,0,s,t,0)$ . And in the first case ,  the following $8$ forms
will suffice :
\begin{equation}
X_{1}^{i}X_{2}^{3-i}X_{4}^{3}(0\leq i\leq 3),
X_{1}^{i}X_{3}^{3-i}X_{5}^{3}(0\leq i\leq 3).\notag
\end{equation}

In the second case , the following $8$ forms will suffice :

\begin{equation}
X_{1}^{i}X_{2}^{3-i}X_{4}^{3}(0\leq i\leq 3),
X_{1}^{i}X_{2}^{3-i}X_{5}^{3}(0\leq i\leq 3).\notag
\end{equation}

The case of $C_{2}$ and $C_{3}$

This follows from a similar argument as we did in   the last case
 .



\begin{thebibliography}{99}



\bibitem{Ba}V. V. Batyrev , {\em On the classification of smooth projective toric varieties }, Tohoku Math.
J. 43 , 569每585 , 1991 .

\bibitem{Cle}H. Clemens , {\em Homological equivalece , modulo algebraic equivalence , is not finitely
generated }, Publ. Math. I.H.E.S. 58 , 19每38 , 1983 .

\bibitem{Cle2}H. Clemens , {\em Double solids }, Adv . in Math . 47 , 107-230 , 1983 .

\bibitem{Co1}D. A. Cox , {\em  The homogeneous coordinate ring of a toric variety }, J. Algebraic Geom. 4
(1): 17-50 , 1995 .

\bibitem{Co2}D. A. Cox , {\em  The functor of a smooth toric variety }, Tohoku Math. J.
 47(2): 251-262, 1995 .


\bibitem{F1}R. Friedman , {\em Simultaneous resolution of threefold double points
 },Math. Ann. 274(4) : 671-689 , 1986 .

\bibitem{F2}R. Friedman , {\em  On threefolds with trivial canonical bundle }, Proc. Symp. Pure. Math. 53 (1991) .

\bibitem{GH}P.S.Green , T.H\"{u}bsch , {\em Connetting moduli spaces of Calabi-Yau
threefolds }, Comm. Math. Phys. 119 : 431每441 , 1988 .

\bibitem{JoKl} T. Johnsen , S. L. Kleiman
, {\em  Rational curves of degree at most 9 on a general quintic
threefold }, Comm. Algebra. 24(8) : 2721-2753 , 1996 .

\bibitem{JoKn}T.Johnsen , A.L.Knutsen , {\em Rational curves in Calabi-Yau threefolds }, Comm. Algebra. 31(8) : 3917-3953 , 2003 .

\bibitem{KS}M. Kreuzer , H. Skarke , {\em Complete classification of reflexive polyhedra in four dimensions }, Adv. Theor. Math. Phys. 4 , no. 6, 1209-1230, 2000.

\bibitem{LT}P.Lu , G.Tian ,  {\em Complex Structures on Connected Sums of $S^{3}\times S^{3}$},
Manifolds and geometry (Pisa , 1993) , 284-293.

\bibitem{Reid}M. Reid , {\em The moduli space of $3-$folds with $K=0$ may neverthless be irreducible},
 Math.Ann. 287 : 329-334 ,1987.

\bibitem{R}M. Rossi , {\em  Geometric Transitions }, J. Geom. Phys. 56(9), 1940-1983 , 2006. 

\bibitem{Ti}G. Tian , {\em Smoothing 3-folds with trivial canonical bundle and ordinary double points} , Essays on Mirror Manifolds (S.T.Yau ed.) , Hongk Kong : International Press , 458-479 , 1992 .


\bibitem{Wa}T.C.T.Wall , {\em  Classification problems in differential topology V. On certain
6-manifolds} , Invent. Math. (1) 355-374 , 1966 ; corrigendum ibid .
(2) 306 , 1967.









\end{thebibliography}
\end{document}